\documentclass[article,10pt,oneside]{memoir}

\setlrmarginsandblock{15mm}{15mm}{*}
\setulmarginsandblock{18mm}{18mm}{*}
\checkandfixthelayout

\usepackage[T1]{fontenc}
\usepackage{amsmath,amssymb}
\usepackage{mathtools}
\usepackage{graphicx}
\usepackage{booktabs}
\usepackage{hyperref}
\usepackage[dvipsnames]{xcolor}
\usepackage{algorithm}
\usepackage{algorithmic}
\hypersetup{colorlinks=true, citecolor=teal, linkcolor=BrickRed, urlcolor=NavyBlue}

\usepackage{fontspec}
\usepackage{unicode-math}
\usepackage{microtype}
\newtoggle{usetikz}

\togglefalse{usetikz}

\iftoggle{usetikz}{
\usepackage{tikz, pgfplots}
\pgfplotsset{compat=1.18}
\usepgfplotslibrary{groupplots}
\usepgfplotslibrary{fillbetween}
\usetikzlibrary{backgrounds}
\usetikzlibrary{arrows.meta}
\pgfplotsset{plot coordinates/math parser=false}
\usetikzlibrary{external}
\usepgfplotslibrary{patchplots}
\usetikzlibrary{shapes.geometric}
\usetikzlibrary{calc}
\usetikzlibrary{backgrounds}
\usetikzlibrary{intersections}
\usetikzlibrary{spy}

\tikzexternalize[mode=list and make, prefix=tikzcache/]
}{}

\renewcommand{\thesection}{\arabic{section}}
\setsecnumdepth{subsubsection}
\numberwithin{equation}{section}

\newcommand{\R}{\mathbb{R}}

\newcommand{\vi}{{\circ}}

\begin{document}

\twocolumn[{%
  \begin{center}
    {\Large\bfseries A one-dimensional model of catheter aspiration in blood vessels\par}
    \vspace{0.5em}
  Niklas Kolbe\textsuperscript{1},
  Dimitrios Katsaounis\textsuperscript{1},
  Xiangge Yin\textsuperscript{2},
  Michael Neidlin\textsuperscript{2}\par
    \vspace{0.3em}
    {\small
      \textsuperscript{1}Institute of Geometry and Applied Mathematics,  RWTH Aachen University, Germany\\
      \textsuperscript{2}Department of Cardiovascular Engineering, RWTH Aachen University, Germany\par}
    \vspace{1em}

    \begin{minipage}{0.85\textwidth}
      \begin{center}
        \bfseries Abstract
      \end{center}

Aspiration thrombectomy is a primary endovascular treatment option for acute ischemic stroke caused by large vessel occlusions, in which a catheter is navigated to the occlusion site and suction is applied to remove the thrombus. This work introduces a reduced one-dimensional model for aspiration in blood vessels that accounts for the elasticity of both the vessel wall and the catheter. The inclusion of vessel wall visco-elasticity transforms the governing equation for the flow rate into a parabolic form, enabling an accurate resolution of the diffuse pressure gradients near the catheter tip observed in CFD simulations. The numerical treatment is based on a relaxation of the hyperbolic subsystem that facilitates nodal solvers, enabling efficient coupling between catheterized and uncatheterized vessel segments, including bifurcations and the catheter tip. A splitting strategy ensures that the visco-elastic terms incur only negligible additional computational cost relative to the purely hyperbolic model. The model is benchmarked against three-dimensional CFD simulations and reference data from the literature. Numerical experiments investigating the role of catheter elasticity and suction force on the hemodynamics are presented, and an uncertainty quantification study demonstrates the suitability of the framework for efficient parameter studies. The results indicate that the proposed model resolves blood flow comparably accurate to three-dimensional simulations at a fraction of the computational cost. While further in vitro and in vivo validation is required, the framework represents a promising step toward computationally efficient preoperative planning and uncertainty-aware decision support in aspiration thrombectomy.

      \medskip
      \noindent\textbf{Keywords:} 1D blood flow modeling; endovascular thrombectomy; finite volume method; nodal solvers
    \end{minipage}
  \end{center}
  \smallskip
  \vspace{1em}
}]

\section{Introduction}
Acute ischemic stroke is one of the leading causes of death and disability worldwide~\cite{donkor2018strokcentur}. Endovascular thrombectomy (EVT) has recently become the standard treatment option for acute ischemic stroke due to large vessel occlusion. During EVT, the cervical vessel of the occluded artery is catheterized with a long and large-lumen access catheter. Smaller catheters are then passed through this access catheter to the intracranial occlusion site and the thrombus is removed mechanically using various techniques.
To date, stent-retriever thrombectomy  is the most common technique, with aspiration thrombectomy with a distal aspiration catheter (so-called ADAPT technique) and variants or combinations of these techniques becoming more common. Common to all mechanical recanalization procedures is that suction is applied to the access catheter to reverse blood flow and prevent thrombus fragments from floating off~\cite{boisseau2020direc}. Regardless of the specific technique, blood flow control during EVT is necessary for the procedure to be effective and safe~\cite{lally2015in,schoenfeld2020effecballoonguide}. The size of the catheter as well as the position of the
tip and the Circle of Willis (CoW) anatomy are important factors for procedural success: a large catheter
allows for more suction, collaterals in turn redirect flow and impair flow control~\cite{rizvi2018inremotaspir,zhao2019collatwillis}. The procedure is highly complex and patient-specific, motivating the use of computational models.  

Computational fluid dy\-na\-mics (CFD) and fluid-structure interaction (FSI) simulations in three space-dimensions have become valuable tools for investigating the hemodynamics of aspiration thrombectomy. CFD models have been employed to evaluate how the aspiration pressure magnitude alters cerebral blood flow, thereby supporting strategies that balance suction efficacy against the preservation of distal perfusion~\cite{luisi2023invescerebhemod}. A hybrid framework coupling 1D and 3D models was used to quantify the effects of catheter tip location and aspiration flow rate, revealing that the mean middle cerebral artery flow is more sensitive to the applied aspiration rate than to the catheter position~\cite{benemerito2023}. Beyond steady aspiration, time-dependent models of recanalization through cyclic aspiration have been developed, enabling the in silico optimization of aspiration waveforms~\cite{good2019devel} and finite element analysis has been used to investigate failure mechanisms in long clots, suggesting that cyclic pressure loading can improve recanalization outcomes compared with continuous aspiration~\cite{patki2022computacute}. FSI simulations that account for vessel wall compliance have been shown to introduce measurable differences in distal pressure and wall shear stress compared with rigid-wall assumptions~\cite{mansilla2025patienspecifmodel}, which underscores the importance of resolving fluid–structure interactions when local hemodynamic accuracy is required. Despite these advances, the translation of 3D CFD and FSI models into clinical decision-making remains limited by their high computational cost, sensitivity to patient-specific boundary conditions, and a lack of comprehensive clinical validation \cite{bai2026mechan}. This has motivated the development of reduced models that resolve essential hemodynamics while achieving the computational efficiency needed for time-critical applications such as treatment planning.

Despite recent advances in the development of surrogate models and data-driven methods, see e.g.,~\cite{eichinger2022surrog, wu2026digit}, one-dimensional blood flow models remain highly relevant in the computational modeling of hemodynamics. Such reduced models are commonly derived from mass conservation and momentum balance, combined with a simplified wall mechanics assumption of static equilibrium in the radial direction of a cylindrical vessel. This yields a first-order hyperbolic system of partial differential equations governing the cross-section area and the flow rate as functions of time and the axial coordinate, that captures wave propagation and pressure--flow relationships along vascular segments at a fraction of the computational cost of 3D simulations. The mathematical properties of these systems, including well-posedness, smoothness of solutions, and numerical treatment, have been extensively studied \cite{canic2003mathem, peiro2009reduc}. To account for the full vasculature, 1D models have been extended to arterial networks through appropriate coupling conditions at bifurcations and junctions~\cite{wang2014}, and the resulting frameworks have been employed in a variety of applications, including the modeling of vascular prosthetics~\cite{canic2002blood}. Within the context of cerebral hemodynamics, 1D models have been widely used to study pulse wave propagation in the CoW, to assess the hemodynamic impact of arterial stenosis and occlusion on global cerebral blood flow, and to evaluate collateral circulation under pathological conditions~\cite{bai2026mechan}. Occlusions within such network models have moreover been considered in \cite{phan2020computmodelclot}. However, the modeling of endovascular therapy in these settings has so far mostly been limited to numerical experiments varying coupling and boundary conditions.

Several studies have investigated fluid flow in catheterized vessels, though predominantly under the assumption of rigid vessel walls and rigid catheters. Early analytical models derived steady-flow resistance relations for annular sections under various rheological assumptions~\cite{wein1970nonnewton, dash1996estim}, with extensions to curved vessels~\cite{jayaraman1995flow} and stenosed arteries with catheters~\cite{back1996flowratepress}. Unsteady analyses include the work of Tsangaris~\cite{tsangaris1984oscil} for a rigid catheter in a rigid vessel and Sarkar and Jayaraman \cite{sarkar2001nonlin}, who considered an elastic vessel with an elastic catheter under small-deformation assumptions. Kanai et al.~\cite{kanai1970one} analyzed wave reflections at the catheter tip for a rigid catheter in an elastic artery, though under simplifying assumptions regarding the phase relationship between pressure and wall velocity. While these analytical contributions provide valuable insight into the mechanics of catheterized flow, they are not directly suited for network-scale simulations. In the context of numerical 1D modeling, Neidlin et al.~\cite{neidlin2015} employed a 1D model to simulate flow in the circle of Willis with an aspiration catheter, representing each artery as a rigid resistance and incorporating the catheter effect by modifying the vessel radius without explicitly resolving the annular cross-section. Pradhan et al.~\cite{pradhan2024} made a significant step forward by deriving a complete 1D hyperbolic formulation for blood flow in an elastic vessel with a rigid catheter, including an analysis of the characteristic system and a Discontinuous Galerkin solver capable of handling catheterized bifurcations and catheter-to-catheter-free transitions. Their model was validated against 3D CFD simulations and demonstrated good agreement in pressure drops and flow rates across idealized and realistic arterial geometries. However, the model is limited to rigid catheters and does not include aspiration effects, which were identified as a direction for future work. 

To overcome current limitations this work introduces a new one-dimensional model for aspiration thrombectomy in an elastic blood vessel using an elastic catheter. The model extends the hyperbolic formulation of Pradhan et al. \cite{pradhan2024} by incorporating vessel wall visco-elasticity, which enables accurate resolution of the localized pressure dynamics at the catheter tip on the millimeter scale. A simplified catheter equilibrium approximation is proposed that reproduces the fully elastic catheter model at a fraction of the computational cost. The numerical treatment relies on a relaxation-based finite volume scheme that, through simple nodal solvers, enables straightforward coupling between catheterized and uncatheterized vessel segments, including bifurcations and the catheter tip. The model is benchmarked against 3D CFD simulations and reference data from the literature in various test cases, including a suction-force--suction-distance analysis and an uncertainty quantification study. 

\section{Methods}
The proposed framework consists of the governing equations derived in Section~\ref{sec:modeling}, the numerical treatment discussed in Sections~\ref{sec:splitting} and \ref{sec:spacedisc} and the coupling of different segments described in Sections~\ref{sec:coupling-problem} and \ref{sec:coupling-problem-catheterized}. In addition, we propose a catheter equilibrium approximation in Section~\ref{sec:catheterequilibrium} and discuss uncertainty quantification by Monte-Carlo estimators in Section~\ref{sec:montecarlo}.   

\subsection{Modeling of blood flow in a catheterized vessel}\label{sec:modeling}
\begin{figure*}
  \centering
  \iftoggle{usetikz}{\input{catheterized.tikz}}{\includegraphics{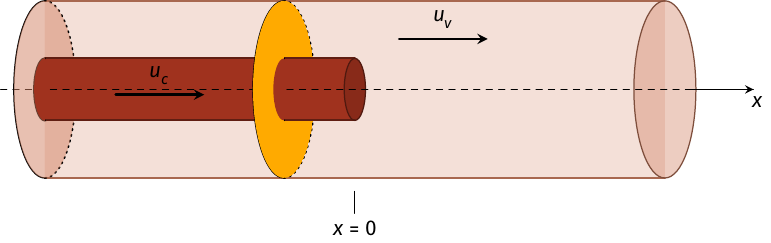}}
  \vspace{2em}
  \iftoggle{usetikz}{\input{cross-section.tikz}}{\includegraphics{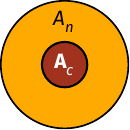}}
\caption{Left: Cylindrical vessel with inserted aspiration catheter ending at position $x=0$. Blood flow is considered both through the catheter and the vessel with different velocities ($u_c$ in the device, $u_n$/$u_v$ in the catheterized/catheter-free vessel). Lower right: Cross section of the catheterized vessel consisting of the round catheter cross-section with area $A_c$ and the annulus-shaped net cross section of the vessel with area $A_n$ subtracting the catheter and its wall from the total cross section.}\label{fig:catheterized}
\end{figure*}

Blood flow is considered through a cylindrical elastic vessel, whose axis is oriented along the $x$-coordinate in the three-dimensional space. Intersecting the vessel with a plane perpendicular to the flow direction at $x=cst$ and neglecting the vessel walls gives rise to the time dependent circular cross section, which we denote as $S(x, t)$. In the center of the vessel an aspiration catheter device is inserted; it is also assumed to be of cylindrical shape and its intersection with the plane is denoted by $\bar S_c(x,t)$. The catheter is itself an elastic pipe with a wall, whose thickness we denote by $h_c$. Consequently, the catheter wall describes an annulus within $S$ and we denote by $S_c(x,t)$ the circle included in $\bar S_c(x,t)$ excluding the catheter wall. Finally we denote by $S_n(x,t)= S(x,t) \setminus \bar S_c(x,t)$ the net vessel cross section, through which the blood of the vessel flows. The corresponding section areas are referred to by $A_n(x,t) \coloneqq |S_n(x,t)|$ and  $A_c(x,t) \coloneqq |S_c(x,t)|$. We refer to Figure~\ref{fig:catheterized} for a visualization of the catheterized configuration.

We assume that blood flows through both, the vessel and the catheter. Following the steps in~\cite{hughes1973} we employ Reynold's transport theorem to derive the catheterized model. Let $\mathbf v$ denote the blood velocity either within the vessel or the catheter then taking the geometry into account the theorem states that for any continuous function $f$ in time and space it holds
\begin{equation}\label{eq:reynold}
\begin{aligned}
  \frac{\partial}{\partial t} (|S_\vi| \, \overline{f\,})+ \frac{\partial}{\partial x} (|S_\vi| \, \overline{f v_1})&= \int_{S_\vi} \left(\frac{\partial f}{\partial t} + \nabla \cdot (f \mathbf v) \right) \, dx \\
  &\quad + \int_{\partial S_\vi} f (\mathbf v_w - \mathbf v) \cdot \mathbf n_\vi \, dS 
\end{aligned}
\end{equation}
for $S_\vi \in \{ S_n, S_c\}$, where $v_1$ refers to the flow velocity in axial direction, $\mathbf n_\vi$ the outer normal of $S_\vi$, $\mathbf v_w$ the velocity of the vessel/device wall and the notation
\[
  \overline{f\,} = \frac{1}{|S_\vi|} \int_{S_\vi} f \, dx   
\]
is used for the section average. We model blood as an incompressible fluid and assume that the vessel and the catheter walls are impermeable and thus $v_w = v$ holds in $\partial S_\vi$. The latter assumption implies that the normal velocities at the device wall from both sides are coupled as $\mathbf v(x, r_c \mathbf n) \cdot \mathbf n = \mathbf v(x, (r_c + h_c) \mathbf n) \cdot \mathbf n$ for any unit vector $\mathbf n$ embedded in the plane $x=cst$ and $r_c\coloneqq \sqrt{A_c/ \pi}$ referring to the local radius of the catheter. Substituting $f=1$ in~\eqref{eq:reynold} now implies the governing equations for the cross section areas of the device and the net vessel, which read
\begin{align}
  \frac{\partial A_\vi}{\partial t} + \frac{\partial A_\vi u_\vi}{\partial x} = 0 \label{eq:At}
\end{align}
with $u_\vi$ referring to the blood velocity in axial direction ($\mathbf v_1$ in \eqref{eq:reynold}) and  $\vi\in\{n, c \}$.

The model is closed by substituting $f=\overline{v_1}$ in \eqref{eq:reynold}, thereby deriving governing equations for the mass flow in the vessel and the catheter defined as $Q_n=A_n u$ and $Q_c=A_c u$, respectively. Key in the derivation, which follows the same path as in the catheter free case, is a reformulation of the material derivative in terms of the Cauchy stress tensor, details are given in~\cite{peiro2009reduc}. We model viscous forces by a simple linear model and assume plug flow, i.e., a flat velocity profile with Coriolis coefficient $\alpha=1$ satisfying
\begin{equation}\label{eq:flatv}
  \overline{v_1^2} = \overline{v_1}^2,
\end{equation}
which provides a suitable approximation for blood flow in large elastic devices and vessels, see~\cite{massey2018mechanfluid}. Moreover, we assume the absence of body forces and obtain the relation
\begin{equation}\label{eq:momentum}
  \frac{\partial}{\partial t} (|S_\vi| \overline{v_1}) + \frac{\partial }{\partial x} \left( |S_\vi|\overline{ v_1}^2 \right) + \frac{|S_\vi|}{\rho} \frac{\partial p_\vi}{ \partial x} + R_\vi \overline {v_1} = 0,
\end{equation}
where $R_\vi$ is a positive viscous resistance parameter and $p_\vi$ is the average pressure in the vessel or the catheter, respectively. Simple algebra allows to deduce from \eqref{eq:momentum} governing equations for the axial velocities that take the form
\begin{align}
  \frac{\partial u_\vi}{\partial t} + u_\vi \frac{\partial u_\vi}{\partial x} + \frac{1}{\rho} \frac{\partial p_\vi}{\partial x} + R_\vi \frac{u_\vi}{A_\vi} = 0. \label{eq:ut}
\end{align}
Unlike the model derivation from the compressible Navier--Stokes equations, see e.g.,~\cite{canic2003mathem}, which makes use of the streamline condition relating axial and radial velocities, here the velocities in different directions remain independent. Consequently the model given by equations~\eqref{eq:At} and \eqref{eq:ut} does not reflect the coupled velocity at the device wall; a coupling between the vessel and the device related quantities is given through the pressure function.

We remark that the flat velocity profile assumed in~\eqref{eq:flatv} implies a lack of spatial variation of the velocity within the cross-section, and makes the model independent of the exact position of the catheter within the catheterized vessel.  

\paragraph*{The uncatheterized vessel} It is instructive to consider the pressure at first in an uncatheterized vessel. In analogy to the above discussion, let $A_v$ denote the section area of a catheter free cylindrical vessel; clearly this section area and the flow velocity in the vessel satisfy the equations~\eqref{eq:At} and~\eqref{eq:ut} taking $\vi=v$. We adopt the independent ring model~\cite{quarteroni2004mathem} for the pressure assuming, among others, linear elastic behavior of the vessel wall and dominance of circumferential stresses. Including the visco-elasticity of the vessel wall in 1D models has improved their accuracy in comparison to in-vitro data in~\cite{alastruey2011pulsewavepropag}. To account for visco-elasticity we employ the Voigt--Kelvin type model~\cite{fung1994biomec} and obtain the relation
\begin{equation}\label{eq:ringmodel}
  \begin{aligned}
   p_v - p_v^\text{ext}  =  \frac{\rho_w h_{v} \,\sqrt{A_{0}}}{\sqrt{A_{v}}} \, \frac{\partial^2 \eta_v}{ \partial t^2} &+ \frac{\gamma_v h_v \pi}{\sqrt{A_v \,A_{0}}} \frac{\partial \eta_v}{\partial t} \\
    &+ \frac{\pi h_v E_v}{(1- \nu^2) \sqrt{A_v \,A_{0}}} \eta_v
  \end{aligned}
\end{equation}
between the pressure $p_v$ within the vessel and the wall displacement $\eta_v$, which due to the cylindrical geometry is given by $\eta_v = \frac{1}{\sqrt{\pi}}(\sqrt{A_v} - \sqrt{A_{0}})$ with $A_{0}$ referring to the cross section area of the vessel in reference configuration. In~\eqref{eq:ringmodel} the parameters $\rho_w$, $\gamma_v$, and $E_v$ refer to the density of the vessel wall material, the viscosity, and the Young modulus of the vessel wall, respectively. In addition, $\nu$ is the Poisson ratio and $p_v^\text{ext}$ the external pressure in the vessel, which we assume constant. Noting that the wall thickness $h_v$ is small in relation to the other quantities in the first term of~\eqref{eq:ringmodel} we proceed as in~\cite{formaggia2003one} and neglect the second order term for simplicity. Another common simplification, see~\cite{quarteroni2004mathem}, is obtained approximating $A_v$ by $A_{0}$ in the linear elasticity term.
Consequently, applying the chain rule we obtain 
\begin{equation}\label{eq:pressure}
  \frac{\gamma_v h_v \sqrt{\pi}}{2 A_v \sqrt{A_{0}}} \frac{\partial A_v}{\partial t} + \beta_v (\sqrt{A_v} - \sqrt{A_{0}}) = p_v - p_v^\text{ext}
\end{equation}
as closed form for the pressure in the vessel, where $\beta_v = \frac{\sqrt{\pi} h_v E_v}{(1 - \nu^2) A_{0}}$. Letting $\tilde p_v=\beta_v (\sqrt{A_v} - \sqrt{A_{0}})$ account for the linear elasticity, the spatial derivative of the pressure is given by
\begin{equation}\label{eq:dpressure}
  \frac{\partial p_v}{\partial x}  =  \frac{\partial \tilde p_v}{\partial x} - \frac{\gamma_v h_v \sqrt{\pi}}{2 A_v \sqrt{A_{v,0}}} \frac{\partial^2 Q_v}{\partial^2 x}
\end{equation}
dropping the small term $ \frac{\partial}{\partial x} \frac{\gamma_v h_v \sqrt{\pi}}{2 A_v \sqrt{A_{0}}}$ and using~\eqref{eq:At}. Taking the time derivative of $Q_v$ and setting $\sigma_v= \frac{\gamma_v h_v \sqrt{\pi}}{2 \rho \sqrt{A_{0}}}$ we end up with the parabolic equation
\begin{equation}\label{eq:Qvt}
  \frac{\partial Q_v}{\partial t} +  \frac{\partial}{\partial x}  \frac{Q_v^2}{A_v}  + \frac{A_v}{\rho} \frac{ \partial \tilde p_v}{\partial x} + R_v \frac{Q_v}{A_v} = \sigma_v \frac{\partial^2 Q_v}{\partial^2 x},
\end{equation}
We note that the equations~\eqref{eq:At} and~\eqref{eq:ut} are equivalent to \eqref{eq:At} and \eqref{eq:Qvt} assuming $\vi=v$, see \cite{sherwin2003one}. We will refer to the former formulation as \emph{velocity form} and to the latter as \emph{flow form} of the blood flow model. If $\sigma_v \ne 0$ the flow form admits a parabolic equation in $Q_v$ for which the numerical solution is straightforward unlike the corresponding velocity form.

Regarding the modeling of the viscous resistance we follow the common literature approach, see e.g., \cite{boileau2015, pimentel-garcia2023high}, and take $R_v=\frac{22 \mu \pi}{\rho}$ with $\mu$ denoting the dynamic blood viscosity, which has been derived from the axial velocity profile proposed in~\cite{smith2002}.

\paragraph*{The catheter} Coming back to the catheterized setting we remark that aspiration catheters are typically made of synthetic polymers like polyurethane or nylon. These polymers exhibit very little visco-elasticity at the molecular level that is not functionally relevant in clinical use~\cite{brydson1999polyet, warner2019guidelearlymanag}. Thus, in the model the pressure in the catheter $p_c$ is only determined by linear elasticity and governed by~\eqref{eq:pressure} taking $v=c$ and $\gamma_c=0$.~i.e.,
\begin{equation}\label{eq:pressurec}
  \beta_c (\sqrt{A_c} - \sqrt{A_{c,0}})   =  p_c - p_c^\text{ext}, \quad \beta_c= \frac{\sqrt{\pi} h_c E_c}{(1 - \nu^2) A_{c,0}}
\end{equation}
with $A_{c, o}$ denoting the reference section area of the catheter. For simplicity viscous resistance in the catheter is modeled in analogy to the vessel taking $R_c= R_v$.  

\paragraph*{The catheterized vessel} The pressure law within the catheterized vessel is also derived from the ring model~\eqref{eq:ringmodel}. Let $A_0$ denote the cross section area of the vessel in reference configuration; note that here the complete vessel including the catheter is considered. The displacement of the vessel wall, which we denote by $\eta_v$, thus satisfies
\begin{equation}\label{eq:etav}
  \sqrt{\pi} \eta_v + \sqrt{A_0} = \sqrt{A_t}
\end{equation}
with
\[
  A_t=|S| = \left(\sqrt{\frac{A_c}{\pi}} + h_c \right)^2 \pi + A_n
\]
referring to the total section area.
The pressure function in the catheterized vessel is based on the displacement~\eqref{eq:etav}. In analogy to the derivation of~\eqref{eq:pressure}, neglecting small temporal changes in $A_c$ as well as the catheter width in the first order term, we obtain from~\eqref{eq:ringmodel} and taking $\beta_v = \frac{\sqrt{\pi} h_v E_v}{(1 - \nu^2) A_0}$ the pressure form
\begin{equation}\label{eq:pressuren}
   p_n - p_n^\text{ext} = \frac{\gamma_n \sqrt{\pi}}{2 A_t \sqrt{A_0}} \frac{\partial A_n}{\partial t} + \beta_v (\sqrt{A_t}  - \sqrt{A_0}).
\end{equation}
Similar to~\eqref{eq:dpressure} and~\eqref{eq:Qvt} this gives rise to a parabolic equation for the net mass flow in the vessel, which defining $\sigma_n= \frac{A_{n,0}}{A_0}\frac{\gamma_v h_v \sqrt{\pi}}{2 \rho \sqrt{A_0}}$ reads
\begin{equation}\label{eq:Qnt}
  \frac{\partial Q_n}{\partial t} +  \frac{\partial}{\partial x}  \frac{Q_n^2}{A_n}  + \frac{A}{\rho} \frac{ \partial \tilde p_n}{\partial x} + R_n \frac{Q_n}{A_n} = \sigma_n \frac{\partial^2 Q_n}{\partial^2 x}.
\end{equation}
after employing a linearization of $A_n$ and $A_c$ and denoting the linear elasticity as $\tilde p_n=\beta_v (\sqrt{A_t} - \sqrt{A_0})$. The reference states of the section area in the catheterized case are related through $A_0=(\sqrt{A_{c,0}/\pi} +h_c)^2 + A_{n,0}$. Motivated by the velocity profile for annular cross sections derived in~\cite{sherwin2003one} we choose the resistance constant
\[
  R_n = \frac{8 \mu \pi (1 - R_0/\sigma_0)}{\rho (2 + (1 - R_0/\sigma_0) ((\log(\sigma_0) - \log(R_0))^{-1} - 1))}
\]
for reference radius $R_0=\sqrt{A_0/\pi}$ and $\sigma_0=\sqrt{A_{c,0} / \pi} + h_c$.

In summary, the full model of the elastic vessel catheterized by an elastic device is given by equation~\eqref{eq:At} for $\circ\in \{n,c\}$ governing the areas $A_n$ and $A_c$, equation \eqref{eq:Qnt} governing the net mass flow $Q_n$ and equation~\eqref{eq:ut} for $\vi=c$ governing the catheter internal blood velocity $u_c$.  

\subsection{Visco-elasticity splitting}\label{sec:splitting}
To handle the second order terms in a way that accounts for the stiffness of the model we employ a splitting approach that we illustrate for the uncatheterized case in the following. We first consider the blood flow model given by~\eqref{eq:At} for $\vi=v$ and~\eqref{eq:Qvt}, neglecting visco-elasticity (taking $\gamma_v=0$) with respect to the variables $A_v$ and $u_v$, which reads
\begin{equation}\label{eq:conservativeAu}\left\{~
  \begin{aligned}
    \frac{\partial A_v}{\partial t} + \frac{\partial A_v u_v}{\partial x} = 0,\\
    \frac{\partial u_v}{\partial t} + \frac{\partial }{\partial x} \left( \frac 12 u_v^2 + \frac{\beta_v}{\rho} (\sqrt{A_v} - \sqrt{A_v^0})\right) + R_v \frac{u_v}{A_v} = 0
  \end{aligned}\right. .
\end{equation}
Visco-elasticity is then accounted for by an isolated heat equation in the variable $Q_v$, i.e.,
\begin{equation}\label{eq:qsplit}
  \frac{\partial Q_v}{\partial t} = \sigma_v \frac{\partial^2 Q_v}{\partial^2 x}.
\end{equation}

\begin{algorithm}[b]
  \caption{Visco-elasticity splitting algorithm for the time discretization of the blood flow model given by \eqref{eq:At} and~\eqref{eq:Qvt}}\label{alg:splitting}
  \begin{algorithmic}[1]
    \REQUIRE initial data $A_v^0$ and $u_v^0$ discretized over a grid
    \STATE Set $t \gets 0$;
    \WHILE{$t < T$}
    \STATE Determine time increment $\Delta t^n$;
    \STATE Compute $A_v^{n+1}$ and $u_v^{n+1/2}$ solving \eqref{eq:conservativeAu} at time $\Delta t^n$ with initial data $A_v^{n}$ and $u_v^{n}$;
    \STATE Compute $Q_v^{n+1}$ solving \eqref{eq:qsplit} at time $\Delta t^n$ with initial data $A_v^{n+1} \, u_v^{n+1/2}$;
    \STATE Set $u_v^{n+1} \gets Q_v^{n+1}/A_v^{n+1}$;
    \STATE Set $n \gets n+1$ and $t \gets t + \Delta t^n$;
    \ENDWHILE
\end{algorithmic}
\end{algorithm}

To compute a numerical solution of the full blood flow system we employ first-order operator splitting alternating between solving the hyperbolic system \eqref{eq:conservativeAu} and the parabolic equation \eqref{eq:qsplit} and taking the solution of~\eqref{eq:conservativeAu} as initial data for~\eqref{eq:qsplit}. Note that variable transformations from $u$ to $Q$ and reverse need to be used in the process. We refer to Algorithm~\ref{alg:splitting} for details. The employed discretization of the hyperbolic system and the parabolic equation is discussed in Section~\ref{sec:spacedisc}.

We employ the same approach to numerically solve the model of the catheterized vessel. In this case we split the model into the hyperbolic system
\begin{equation}\label{eq:conservativecatheter}\left\{~
  \begin{aligned}
    \frac{\partial A_n}{\partial t} + \frac{\partial A_n \, u_n}{\partial x} = 0, \qquad \frac{\partial A_c}{\partial t} + \frac{\partial A_c \, u_c}{\partial x} &= 0, \\
    \frac{\partial u_n}{\partial t} + \frac{\partial }{\partial x} \left( \frac 12 u_n^2 + \frac{\beta_n}{\rho} (\sqrt{A_t} - \sqrt{A_0})\right) + R_n \frac{u_n}{A_n} &= 0, \\
    \frac{\partial u_c}{\partial t} + \frac{\partial }{\partial x} \left( \frac 12 u_c^2 + \frac{\beta_c}{\rho} (\sqrt{A_c} - \sqrt{A_{c,0}})\right) + R_c \frac{u_c}{A_c} &= 0
  \end{aligned}\right.
\end{equation}
in the variables $A_n$, $u_n$, $A_c$, and $u_c$ and a heat equation in $Q_n$ similar to \eqref{eq:qsplit} and proceed in analogy to Algorithm~\ref{alg:splitting}. The main motivation for employing the velocity form within this splitting approach is the handling of the coupling problem that we expand upon in Section~\ref{sec:coupling-problem}.

\subsection{Numerical discretization}\label{sec:spacedisc}
We use the approach developed in~\cite{beckers2025laxfried} to solve the hyperbolic systems~\eqref{eq:conservativeAu} and~\eqref{eq:conservativecatheter}. In short, given a fixed discretization of the vessel into computational cells and a time increment $\Delta t$ the Lax--Friedrichs-type scheme
\begin{equation}\label{eq:lxf}
  \mathbf U^{n+1}_j = \mathbf U^n_j - \frac{\Delta t}{\Delta x}\left(\mathcal{F}_{j+1/2}^n - \mathcal{F}_{j-1/2}^n \right) - \Delta t\, \mathbf S(\mathbf U_j^n)
\end{equation}
with numerical fluxes
\[
  \mathcal{F}_{j-1/2}^n = \frac 12 (\mathbf V_{j-1}^n + \mathbf V_j^n) - \frac{1}{2} \Lambda (\mathbf U^n_{j} - \mathbf U_{j-1}^n)
  \]
updates the cell average of the system state $\mathbf U$ over cell $j$ from time $t^n$ to time $t^{n+1}=t^n+ \Delta t$. Away from the couplings and domain boundaries the variable $\mathbf V$ coincides with the flux function evaluated at the system state, i.e. $\mathbf V_j^n = \mathbf F(\mathbf U_j^n)$.
System state, flux function and source term either take the form
\begin{align*}
  \mathbf U = (A_v, u_v)^T, \\
  \mathbf F(\mathbf U) = \left(A_v\, u_v, \frac 12 u_v^2 + \frac{\beta_v}{\rho} (\sqrt{A_v} - \sqrt{A_{v,0}})\right)^T, \\
  \mathbf S(\mathbf U) = \left(0, K_R \frac{u_v}{A_v}\right)^T
\end{align*}
if system~\eqref{eq:conservativeAu} is addressed or
\begin{align*}
  \mathbf U &= (A_n, u_n, A_c, u_c)^T, \\
  \mathbf F(\mathbf U) &= \Big(A_n\, u_n, \frac 12 u_n^2 + \frac{\beta_n}{\rho} (\sqrt{A_t} - \sqrt{A_0}),\\
  &\qquad A_c\, u_c, \frac 12 u_c^2 + \frac{\beta_c}{\rho} (\sqrt{A_c} - \sqrt{A_{c,0}})\Big)^T, \\
  \mathbf S(\mathbf U) &= \left(0, K_R \frac{u_n}{A_n}, 0, K_R \frac{u_c}{A_c} \right)^T
\end{align*}
if system \eqref{eq:conservativecatheter} is addressed. The Matrix $\Lambda$ is constant over all cells and determined by the subcharacteristic condition of the system, see~\cite{chen1994hyperconserlaws}. In practice we take $\Lambda = \lambda I$ in the uncatheterized case with $\lambda$ referring to the maximal characteristic velocity, which is updated in every time step together with the time increment $\Delta t \propto \lambda^{-1}$, see~\cite{beckers2025laxfried} for details. In the catheterized setting we extend this approach and maximize two separate characteristic velocities: the one in the vessel (governing $\lambda$) and the one in the catheter (governing $\lambda_c$). The diagonal matrix $\Lambda$ then uses the entry $\lambda$ in entries that correspond to vessel-related components 
and $\lambda_c$ in entries that correspond to catheter-related components. In this way, we account for the different magnitude of the two speeds and significantly reduce numerical viscosity over a unified treatment of the model equations.

The parabolic equation \eqref{eq:qsplit} is numerically solved on the same grid as the hyperbolic system using the standard three-point finite difference formula combined with the backward Euler method. The overall scheme is of first order accuracy, which is sufficient to accurately resolve the test cases on fine meshes. Although this is not considered for the purpose of this work, we note that a scheme of second order accuracy can be constructed using Strang splitting instead of Godunov splitting in Algorithm~\ref{alg:splitting} and enhancing~\eqref{eq:lxf} using piecewise linear reconstructions, see~\cite{lanser1999analy, leer1979towarultimconserdifferschem}. 

Due to the convection dominant nature of the problem, we employ characteristics-based boundary conditions derived from the non-diffusive systems~\eqref{eq:conservativeAu} and~\eqref{eq:conservativecatheter}, that assume no visco-elasticity at the boundaries. With respect to the splitting proposed in Section~\ref{sec:splitting} this means that homogeneous Neumann conditions are assumed, when the heat equation~\eqref{eq:qsplit} is solved, while experiment specific conditions for the hyperbolic systems are used. This boundary handling has been commonly used in computational fluid dynamics, see, e.g.,~\cite{wesseling2001princ}. We use in particular resistance conditions with large reflection coefficients to account for the blockage by a thrombus in the numerical experiments. Inlet pressure and flow conditions in the catheterized vessel are realized in two steps computing boundary data first for the isolated catheter and employing this information to compute the boundary data regarding the variables $A_n$ and $u_n$ following~\cite{beckers2025laxfried}.

\subsection{Coupling at the tip of the device}\label{sec:coupling-problem}
In this section we consider the transition between catheterized and un-catheterized segments within a vessel. We assume that a vessel parameterized by its axial position $x \in \R$ has the aspiration device inserted on the negative half axis $(-\infty, 0)$ affecting also the blood flow in the catheter free segment $(0, \infty)$.  We refer again to Figure~\ref{fig:catheterized} for a visualization of the configuration.

We first consider the visco-elasticity free case ($\gamma_n=0$), in which this situation is modeled coupling model~\eqref{eq:conservativecatheter} on the negative half-axis with~\eqref{eq:conservativeAu} on the positive half-axis. For closure we impose continuity of the total cross section radius, i.e.
\begin{equation}\label{eq:couplingA}
  \sqrt{A_t(0^-, t)} = \sqrt{A_v(0^+, t)},
\end{equation}
as first coupling condition and note that $A_t$ depends on the areas $A_n$ and $A_c$ through \eqref{eq:etav}. As second coupling conditions we impose the continuity of the total mass flow, i.e.,
\begin{equation}\label{eq:couplingu}
  A_c(0^-, t) \, u_c(0^-, t) + A_n(0^-, t) \,u_n(0^-, t) = A_v(0^-, t)\,u(0^+, t).
\end{equation}
In addition, we assume homogeneous boundary conditions regarding the catheter related quantities at the interface, which are given by
\begin{equation}\label{eq:couplinghomogeneous}
  \frac{\partial A_c(0^-, t)}{\partial x} =  0 \quad \text{and} \quad \frac{\partial u_c(0^-, t)}{\partial x} = 0.
\end{equation}

To solve this coupling problem we rely on the framework in~\cite{herty2023centrschemtwo} and consider Jin-Xin relaxation of the two systems, see~\cite{jin1995relaxschemsystem}, giving rise to additional relaxation variables. After identifying consistent coupling conditions for those variables two coupled half-Riemann problems are solved explicitly delivering Dirichlet boundary data for scheme~\eqref{eq:lxf} at the interface. We refer to Appendix~\ref{sec:hcoupling} for further details.

If visco-elasticity is assumed we follow a hybrid approach employing the splitting proposed in Section~\ref{sec:splitting}. In the hyperbolic step analogue to line~4 of Algorithm~\ref{alg:splitting} we follow the above procedure, which particularly identifies the coupling states $[Q_n]_R= [A_n]_R \, [u_n]_R$ and $[Q_v]_L = [A_v]_L \, [u_v]_L$, cf. Appendix~\ref{sec:hcoupling}. In the diffusion step analogue to line~5 of Algorithm~\ref{alg:splitting} we decompose the catheterized and the un-catheterized subdomain and impose Neumann boundary conditions at the interface with slope $([Q_v]_L - [Q_n]_R) / \delta $ implying a transition region of length $\delta$, which we choose as $\delta = 1 $mm in our experiments.
\begin{figure*}
  \iftoggle{usetikz}{\input{1-1.tikz}}{\includegraphics{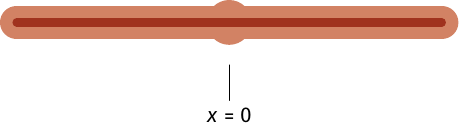}}
  \hspace{1cm}
  \iftoggle{usetikz}{\input{1-2.tikz}}{\includegraphics{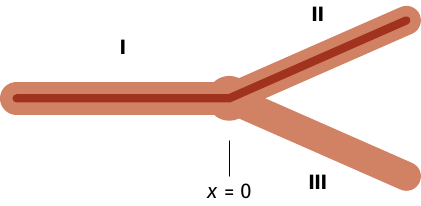}}
\caption{Left: Coupling of two catheterized vessel segments. The model accounts for the coupling by imposing algebraic condition at the coupling interface located in $x=0$. Right: Bifurcation in a catheterized vascular network. At the coupling node in $x=0$ two catheterized (I \& II) and one uncatheterized (III) vessel segments meet.}\label{fig:couplings}
\end{figure*}

\subsection{Coupling of catheterized segments}\label{sec:coupling-problem-catheterized}

To model catheterization in vascular networks we need to account for other interactions of catheterized vessel segments apart from the situation at the tip of the catheter discussed in the previous section. We discuss these models for the inviscid model, in which $\gamma_v=0$ and note that in the viscous case similar splitting techniques can be used as in Section~\ref{sec:coupling-problem}. One case of interest is the coupling of two single vessel segments as shown in Fig.~\ref{fig:couplings} (left). Here, two models of the form~\eqref{eq:conservativecatheter} are coupled at the origin, which allows for a discontinuity in the mechanical properties of the vessel or its geometry, i.e.~in the parameters $A_{0}$ or $\beta_v$. While less relevant in practice, this also allows for the modeling of abrupt changes in the catheter material and geometry. Introducing the total pressure $p_{\vi, t} = \frac{\rho}{2} u_\vi^2 + p_\vi$ we follow the network modeling in \cite{quarteroni2004mathem} and impose continuity of mass flow and total pressure but in both the catheter and the vessel, i.e.
\begin{align}
  p_{c,t}(0^-, t) &= p_{c,t}(0^+, t), \quad &Q_{c}(0^-, t) &= Q_c(0^+, t), \label{eq:coupling11c}\\
  p_{n,t}(0^-, t) &= p_{n,t}(0^+, t), \quad &Q_{n}(0^-, t) &= Q_n(0^+, t).\label{eq:coupling11n}
\end{align}
To solve the coupling problem we rely on the introduced relaxation approach. Thereby we can decouple the catheter and the vessel and first compute the catheter related coupling states from \eqref{eq:coupling11c} and employ those in the subsequent computation of the vessel related coupling states. Both steps here follow the procedure discussed in \cite{beckers2025laxfried}, which is similar to Appendix~\ref{sec:hcoupling}, and requires the solution of a nonlinear scalar equation using Newton's method.

Another case of interest is given by catheterized bifurcations as shown in Fig.~\ref{fig:couplings} (right), in which two catheterized and one uncatheterized vessel segments are coupled. We thus use model \eqref{eq:conservativecatheter} to account for the incoming catheterized vessel segment~\textrm{I} and the outgoing catheterized vessel segment~\textrm{II} and model~\eqref{eq:conservativeAu} to account for the catheter free outgoing segment~\textrm{III}. Combining the common network model with \eqref{eq:coupling11c} and using superscripts to indicate the vessel segments we obtain the coupling conditions
\begin{align}
\left\{\begin{aligned}
  p_{c,t}^\mathrm{I}(0^-, t) &= p_{c,t}^\mathrm{II}(0^+, t), \\
  Q_{c}^\mathrm{I}(0^-, t) &= Q_c^\mathrm{II}(0^+, t), 
\end{aligned}\right. \label{eq:coupling12c}\\
\left\{ \begin{aligned}
  p_{n,t}^\mathrm{I}(0^-, t) &= p_{n,t}^\mathrm{II}(0^+, t) = p_{v,t}^\mathrm{III}(0^+, t),  \\
  Q_{n}^\mathrm{I}(0^-, t) &= Q_n^\mathrm{II}(0^+, t) + Q_n^\mathrm{III}(0^+, t).
  \end{aligned}
  \right.\label{eq:coupling12v}
\end{align}
This coupling problem is adressed by the relaxation framework, where again the catheter related subproblem~\eqref{eq:coupling12c} can be decoupled. The vessel related subproblem~\eqref{eq:coupling12v} can be rearranged in terms of a nonlinear system that can be efficiently solved combining the numerical solution of two scalar problems by Newton's method.

\subsection{Catheter equilibrium approximation}\label{sec:catheterequilibrium}
In our experiments we work with inlet pressure condition both for the aspiration catheter and the vessel. In case of the former this allows for a simplified approximation of the flow in the catheter, which due to large characteristic velocities significantly alleviates the computational costs of numerical simulations. Particularly in case of a constant inlet pressure the flow profile attains an equilibrium after a very short time period. Clearly, this equilibrium is uniform in space. Let $p_c$ denote the constant inlet pressure at the aspiration catheter. In this case the equilibrium of the section area corresponds to the inverse pressure, i.e.,
\begin{equation}\label{eq:Acequilibrium} 
  A_c \equiv \left( \frac{p_c - p_c^\text{ext}}{\beta_c} + \sqrt{A_{c,0}} \right)^2.
\end{equation}
Assuming zero flow in the device before an inlet pressure is applied non-reflecting conditions in~\cite{beckers2025laxfried} motivate the axial velocity
\begin{equation}\label{eq:ucequilibrium} 
u_c \equiv  A_c^{-3/4} \sqrt{\frac{\beta_c}{2 \rho}} (A_c - A_{c,0}).
\end{equation}
While the exact steady state of the catheter sub-model under constant inlet pressure conditions depends on the initial section area and the geometry of the device, fixing the section area~\eqref{eq:Acequilibrium} and the velocity~\eqref{eq:ucequilibrium} within the full coupled problem discussed in Section~\ref{sec:coupling-problem} has lead to an accurate approximation of the flow dynamics throughout our numerical experiments. We further note that the resulting simplified catheterized model is identical to the rigid catheter form in~\cite{pradhan2024} if aspiration is neglected.

\subsection{Uncertainty quantification}\label{sec:montecarlo}
Hemodynamics exhibit considerable heterogeneity, not only across patients but also within distinct regions of the vasculature. Consequently, even with advances in medical imaging, substantial uncertainty remains in the choice of model parameters. Here, we investigate the role of parameter uncertainty using a Monte Carlo estimator. Let $\xi$ be the random variable that represents the model uncertainty taking values in a subset of $\R$. The quantity of interest $\mathcal{Q}$ maps a numerical simulation, for which we fix the setup and the spatial and temporal discretization parameters, to a real valued observable. Through the uncertainty $\mathcal{Q}$ can be understood as a function of $\xi$. Taking now $M$ samples of the random variable $\xi_1, \dots, \xi_M$ the Monte Carlo estimator for the expected value of $\mathcal Q$ is
\[
  \frac{1}{M} \sum_{j=1}^M \mathcal{Q}(\xi_j) \approx \mathbb{E}[\mathcal{Q}].
\]
By arranging the observables $\mathcal{Q}(\xi_j)$ for $j=1, \dots, M$ in an increasing series we further obtain the quantile corresponding to probability $q \in [0,1]$ taking the entry $\lfloor M q \rfloor$ of that series.

\section{Results}
The new model is tested in three cases that investigate the dynamics of aspiration within a partially catheterized arterial segment, near a vascular bifurcation and in case of an arterial blockage. The first test case shows the feasibility of the model dynamics in different time scales and verifies the convergence to the catheter equilibrium approximation, the second and third test cases benchmark the model against literature results and 3D CFD simulations. The finite volume scheme that we use in combination with the relaxation based coupling and boundary treatment was verified in the catheter-free case in~\cite{beckers2025laxfried}. Additional information on the numerical discretization in the presented results is provided in Appendix~\ref{sec:numsetup}. 

In all experiments we choose the initial section area of the vessels and the catheter identical to the corresponding reference configuration $A_{\vi,0}$ and set the velocities $u_\vi\equiv 0$. Simulation dynamics are triggered by the choice of coupling and boundary conditions. If not otherwise mentioned, we conduct our simulations assuming blood density $\rho=1.06$ g/ml, Poisson ratio $\nu=0.5$, vessel wall viscosity $\gamma_v=0$ Pa \,s and dynamic blood viscosity $\mu=0.04$ Pa\,s.

\subsection{Aspiration in a partially catheterized arterial segment}\label{sec:partiallycatheterized}
In the first experiments we investigate the effect of suction through a catheter in a catheterized arterial segment. The medium size vessel segment we consider has a length of $12$ cm, a uniform diameter of $0.4$ cm and a wall width of $h_v=0.05$ cm. A catheter with uniform diameter $0.1$ cm and wall width $h_c=0.05$ cm is inserted from the proximal boundary located at position $x=-6$ cm and ranges to the center of the considered domain ($x=0$ cm). Young's modulus of the vessel and the catheter are chosen $E_v=600$ kPa and $E_c=30$ MPa noting that the magnitude of the latter fits to the results of an indentation test with a relevant catheter. External pressure is taken $p_v^\text{ext}= p_n^\text{ext}=4.3$~kPa and $p_c^\text{ext}=490$~kPa. At the distal boundary ($x=6$) we impose fully absorbent boundary conditions, in which the back propagating Riemann invariant is set to zero, see~\cite{sherwin2003comput}.

\paragraph*{Short time scale dynamics}
\begin{figure*}
\iftoggle{usetikz}{\input{microscale.tikz}}{\includegraphics{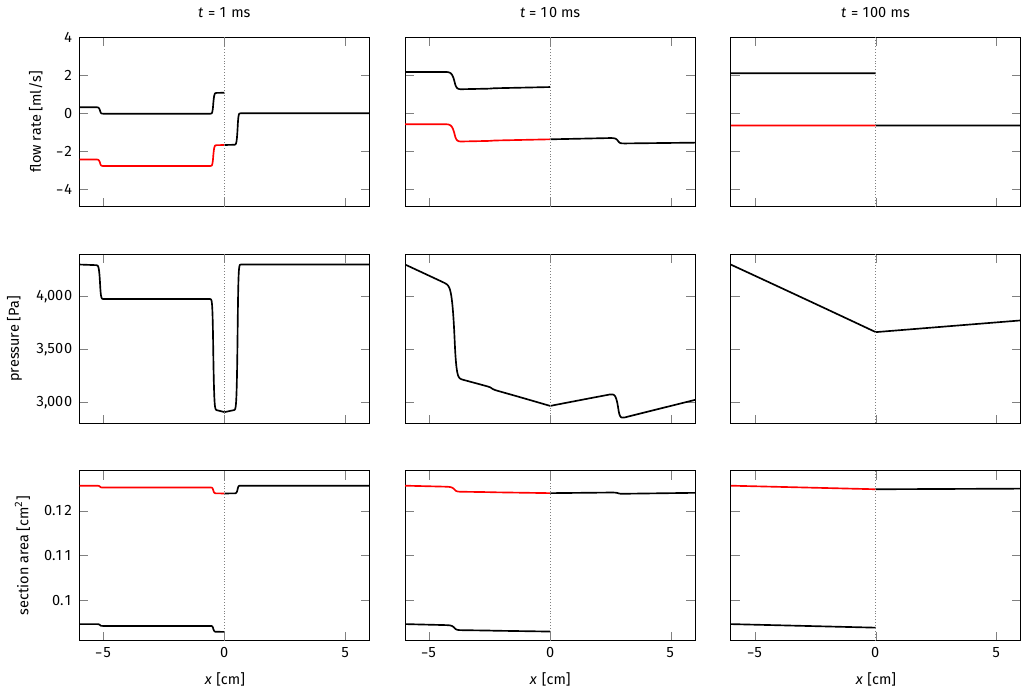}}
\caption{Flow rate, pressure, and section area over a vessel segment at three time instances (black lines) as suction is applied through an aspiration catheter inserted from the left inlet with tip at the center of the segment. Red lines indicate the total flow rate in the top row and the total section area in the bottom row (including the catheter). In the microsecond scale a local decrease of pressure and section area emerges at the tip of the catheter.}\label{fig:microscale}
\end{figure*}
At the left boundary we fix the pressure $p_n^\text{ext}$ at the vessel inlet and apply suction in the catheter imposing the constant pressure $-30$~kPa at the catheter inlet. Figure~\ref{fig:microscale} shows the flow rate, pressure and section area over the vessel segment at three time instances shortly after onset of aspiration. As the negative flow reaches the end of the catheter both the pressure and the section area decrease locally at this position ($t=1$ ms). Due to the elasticity of the vessel the pressure drop is propagated over the full vessel ($t=10$ ms) and eventually a steady state is reached, in which the total flow rate is zero throughout the segment and the pressure attains a convex profile that is linear in the catheterized and the uncatheterized parts of the segment with minimum at the tip location ($t=100$ ms). These short time dynamics, particularly the produced pressure drop, qualitatively match the findings from 3D CFD simulations~\cite{chitsaz2018threedimennumer}.

\paragraph*{Role of the catheter elasticity}
\begin{figure*}
  \iftoggle{usetikz}{\input{elasticity-study.tikz}}{\includegraphics{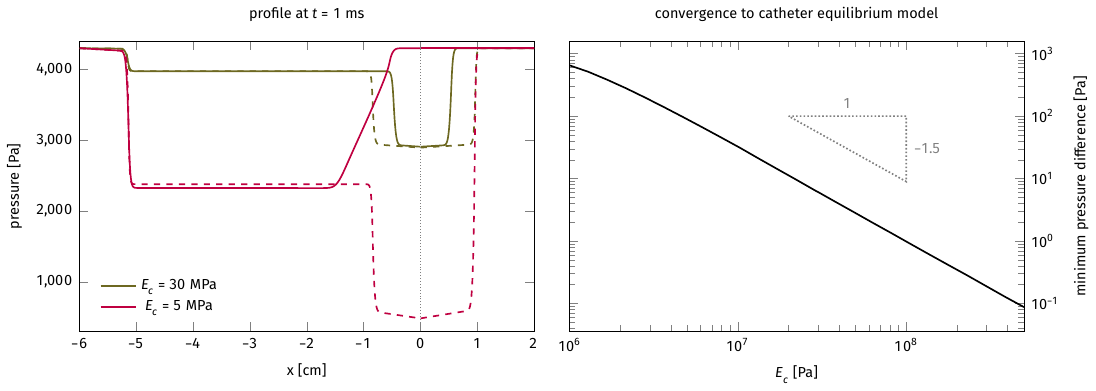}}
  \caption{Left: Comparison of spatial pressure profiles at time $t=1$ ms obtained by the fully elastic catheter model (solid lines) and the catheter equilibrium approximation (dashed line). Two different elastic catheter moduli are considered. Right: Absolute difference of the minimal pressure of the steady state between the fully elastic catheter and the catheter equilibrium model in relation to Young's modulus of the catheter in loglog scale. As the elasticity of the catheter decreases and $E_c$ increases the fully elastic model dynamics converge to the ones of the catheter equilibrium approximation.}\label{fig:elasticity-study}
\end{figure*}

In the same configuration we analyze the behavior of the catheter equilibrium model introduced in Section~\ref{sec:catheterequilibrium} and compare it to the standard model, which we denote in the following as fully elastic model. Numerical computations show that while the equilibrium approximation leads to differences to the fully elastic model in the short time scale (approximately $1-30$ ms) it attains similar steady states on a longer time scale (approximately 100 ms). The spatial pressure profile at $t=1$ ms in the above setup ($E_c=30$ MPa) is compared to the corresponding catheter equilibrium approximation in Figure~\ref{fig:elasticity-study} (left). Although both profiles are of similar magnitude the catheter equilibrium approximation predicts a wider decreased pressure region at the catheter tip in comparison to the fully elastic model. Along Figure~\ref{fig:elasticity-study} (left) presents pressure profiles at the same time for the fully elastic model and the catheter equilibrium approximation in case of the reduced Young's modulus $E_c= 5$ MPa. In this case a significantly larger pressure gradient emerges and the difference between the fully elastic model and the catheter equilibrium approximation are more evident; while the catheter equilibrium approximation leads to a profile that resembles a scaled version of the profiles corresponding to the lower elasticity case the fully elastic model predicts a dissimilar profile as the negative flow wave has not reached the tip position at $t=1$ ms in this model.

For a broader picture we compute simulations using both models for a range of different values for $E_c$. In all computations a steady state was reached before the final time $T=100$ ms similar to the one shown in Figure~\ref{fig:microscale} at the final time. To compare those steady states we compute for fixed $E_c$ the absolute difference of the minimal pressure in the two models, i.e.
\[
  \big| p_\text{eq}(t=100 \text{ ms}, ~x=0 \text{ cm}) - p_\text{el}(t=100 \text{ ms}, ~x=0 \text{ cm}) \big|.
\]
This pressure difference decreases as Young's modulus increases as we show in Figure~\ref{fig:elasticity-study} (right), which also suggest that the pressure difference is proportional to $E_c^{-1.5}$. Altogether, convergence between the two models in the stiff catheter limit is indicated and the equilibrium catheter model appears to be a good approximation of the fully elastic model if Young's modulus of the catheter is high and the dynamics of interest do not occur on a time scale in the order of few microseconds.    

\paragraph*{Medium time scale dynamics}
\begin{figure*}
  \iftoggle{usetikz}{\input{mediumscale.tikz}}{\includegraphics{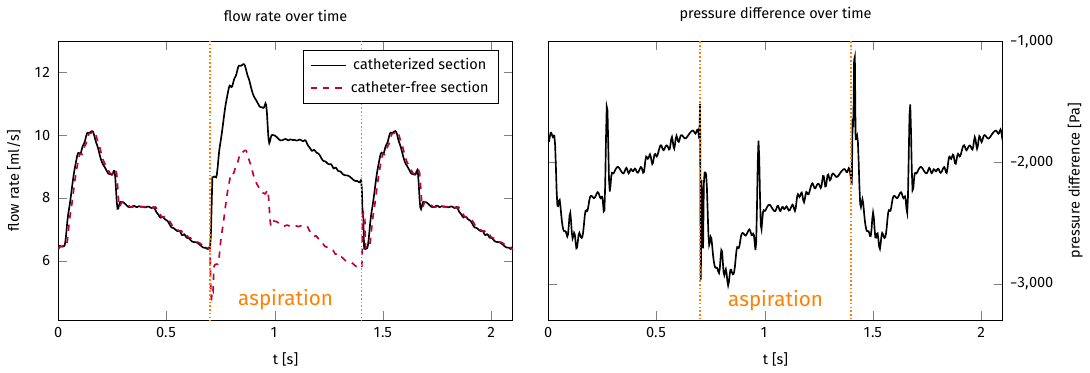}}
  \caption{Flow rate and pressure difference over a time span of three heartbeats with aspiration performed while $t\in[0.7, 1.4]$ (indicated in orange). Left: Flow rate in the vessel in the catheterized (evaluated at $x=-4$ cm) and the catheter-free part of the segment (evaluated at $x=4$ cm).  Right: Pressure difference between the catheter-free (evaluated at $x=4$ cm) and the catheterized part of the segment (evaluated at $x=-4$ cm).}\label{fig:mediumscale}
\end{figure*}

To study the model predictions on a larger time scale we next simulate a period of multiple heart beats by imposing experimental pressure data measured during balloon occlusion of the left internal carotid artery at the vessel inlet. The data comes from an in-vitro model using an elastic phantom of a cohort averaged CoW anatomy~\cite{luisi2024evaluat}. Aspiration is performed during a part of the simulated time window, which is implemented by imposing at the catheter inlet the constant pressure $-30$ kPa during the second heartbeat (as $t\in[0.7, 1.4]$) and $p_c^\text{ext}$ during the rest of the simulated time span. For the calibration of the model we simulate an aspiration-free heartbeat before the relevant time window~$[0, 2.1]$.

Figure~\ref{fig:mediumscale} shows the time evolution of the flow rate in the vessel, i.e., $Q_n$ in the catheterized and $Q_v$ in the catheter-free part of the segment, evaluated at the positions $x=-4$ cm and $x=4$ cm. While the flow rate is nearly constant over the full segment throughout most of the simulated time span, the flow rate in the catheterized part of the segment significantly increases during aspiration and the flow rate in the catheter-free part decreases to a smaller extent. Figure~\ref{fig:mediumscale} moreover shows the pressure difference between the catheter free and the catherized part of the segment over time. This difference is computed again evaluating the pressure in the vessel at $x=4$ cm and $x=-4$ cm. Throughout the considered time span the pressure drops at the transition from the catheterized to the catheter-free part. This drop becomes more pronounced during aspiration; only after offset of aspiration a short reduction of the difference occurs.

Both time curves indicate that blood flow normalizes shortly after offset of aspiration. As an occlusion free scenario is simulated here, this as well as the predicted variation in flow and pressure drop align with observed flow dynamics during aspiration therapy, compare e.g.~\cite{marconi2021invespressflow}.  

\subsection{Aspiration in a catheterized arterial bifurcation}\label{sec:bifurcationstudy}
\begin{figure*}
\begin{tabular}{cc}
\centering
\resizebox{0.45\linewidth}{!}{\iftoggle{usetikz}{\input{bif1.tikz}}{\includegraphics{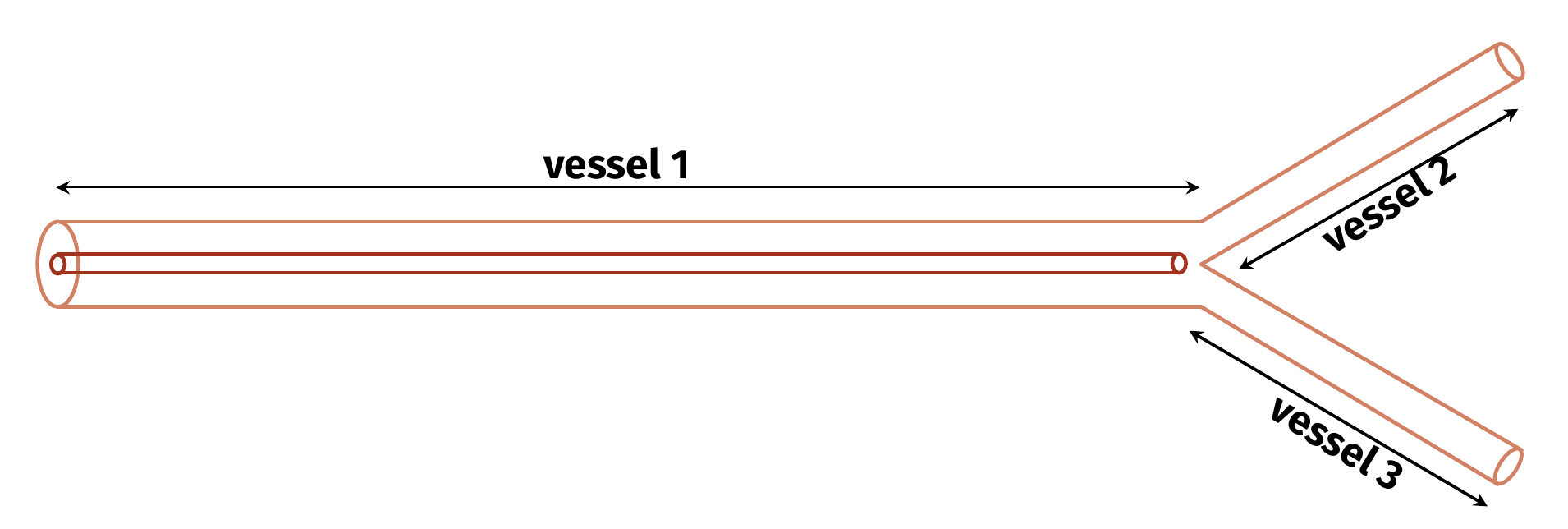}}} & \resizebox{0.45\linewidth}{!}{\iftoggle{usetikz}{\input{bif2.tikz}}{\includegraphics{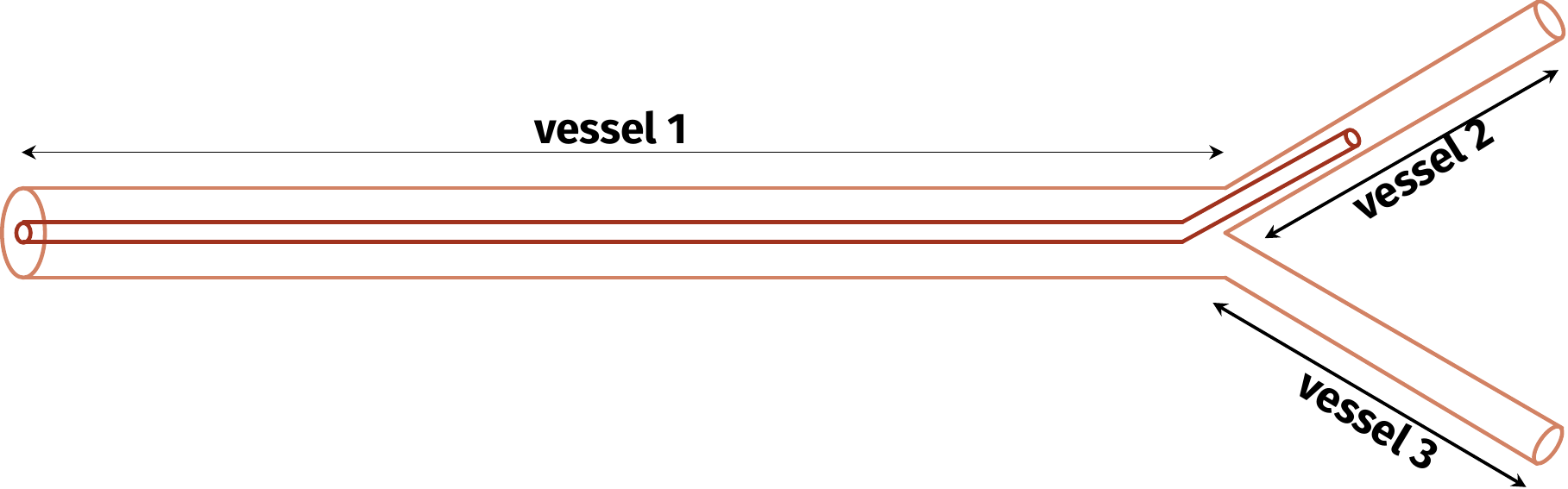}}} \\
  (A) & (B)
\end{tabular}
\caption{Two experimental cases for aspiration in a catheterized arterial bifurcation. (A) The catheter is inserted from the inlet boundary of vessel~1 and stops at the bifurcation point. (B) The catheter is inserted from the inlet boundary of vessel~1 and passes through the bifurcation point into vessel~2.}
\label{fig:bifurcation_schematic}
\end{figure*}

The aim of the next experiment is to test the effects of aspiration on small catheterized arterial networks. In the test case a single arterial vessel (\textit{vessel~1} with inner diameter $0.24$ cm and total length $16$ cm) bifurcates into two daughter vessels (both with length $5$ cm) of different diameter ($0.22$ cm in \textit{vessel~2} and $0.18$ cm in \textit{vessel~3}), as shown in Figure~\ref{fig:bifurcation_schematic}.
The thickness of the arterial walls amounts to $5\%$ of the respective inner vessel diameter and Young's modulus of $E=600$ kPa is assumed for all vessels. 
A catheter (inner diameter $0.12$ cm and wall thickness $0.006$ cm) is inserted from the inlet of the vessel, which either stops directly in front of the bifurcation point (see Figure~\ref{fig:bifurcation_schematic}A), or passes through the bifurcation into one of the two branching vessels, where it eventually ends (see Figure~\ref{fig:bifurcation_schematic}B).
The basic setup was adopted from~\cite{pradhan2024}, where the middle cerebral artery bifurcation was modeled, and a 1D blood flow model with a rigid catheter was validated by a 3D CFD solver in various experiments evaluating the effects of catheterisation on the overall arterial flow.
Here we compare the simulations of our modeling approach with the results by the authors and in addition, investigate the effect of aspiration. 

In the configuration an input flow rate of $0.475\text{ ml/s}$ at the inlet of \textit{vessel 1} and resistence boundary conditions at the outlet of \textit{vessel 2} and \textit{vessel 3} are imposed with resistence $46.3\text{ kPa\,s/ml}$ and $84.1\text{ kPa\,s/(ml)}$, respectively. For consistency, blood density $\rho = 1 \text{ g/ml}$ and no external pressure in the vessels is assumed.
In the first test case, shown in Figure~\ref{fig:bifurcation_schematic}A, the catheter length corresponds with the length of vessel~1 and in the second test case, shown in Figure~\ref{fig:bifurcation_schematic}B, the catheter extends for an additional $2\text{ cm}$ into vessel~2.

Our model assumes an elastic catheter with Young's modulus taken $E_c = 60 $ MPa. We consider two aspiration scenarios, in which a negative inflow rate is imposed at the catheter inlet, that is either taken $-0.05\text{ ml/s}$ (\emph{weak} aspiration) or $-0.25\text{ ml/s}$ (\emph{strong} aspiration).

\begin{figure*}
  \iftoggle{usetikz}{\input{bifurcation_1.tikz}}{\includegraphics{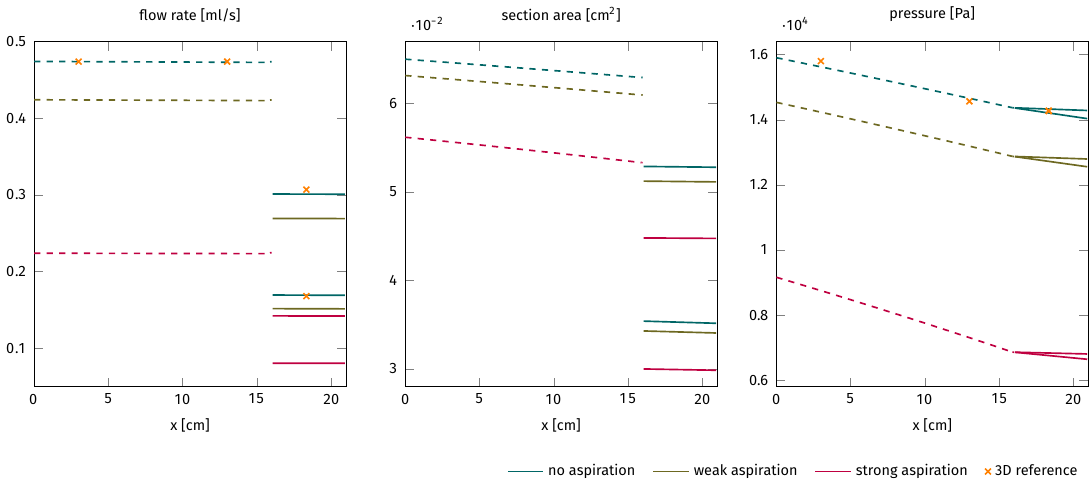}}
  \caption{Model results for the test case in Figure~\ref{fig:bifurcation_schematic}(A), where the catheter stops at the bifurcation point. Flow rate, section area, and pressure over the vessels are presented. In the catheterized segments (represented by dashed lines) total flow rate and section area are shown that include the flow and the size of the catheter. Either no aspiration, weak aspiration (suction with flow rate $-0.05\text{ ml/sec}$) and strong aspiration (suction with flow rate $-0.25\text{ ml/sec}$) is applied.}\label{fig:aspiration_bifurcation}
\end{figure*}

In Figure~\ref{fig:aspiration_bifurcation}, we present simulation results of the setting in Figure~\ref{fig:bifurcation_schematic}A, showing total flow rate, section area and pressure, over the three vessels, when an equilibrium state is reached.
In the aspiration free cases, the model agrees with the 3D CFD literature results, marked in Figure~\ref{fig:aspiration_bifurcation} as well as with the 1D model proposed by the authors. 
The application of aspiration therapy leads to a significant flow rate reduction in the arterial vessels, as well as a decrease of the section area throughout the catheterized segment. The effects of both catheterization and aspiration are more pronounced closer to the tip of the catheter, where the vessel shrinks significantly.  

\begin{figure*}
  \iftoggle{usetikz}{\input{bifurcation_2.tikz}}{\includegraphics{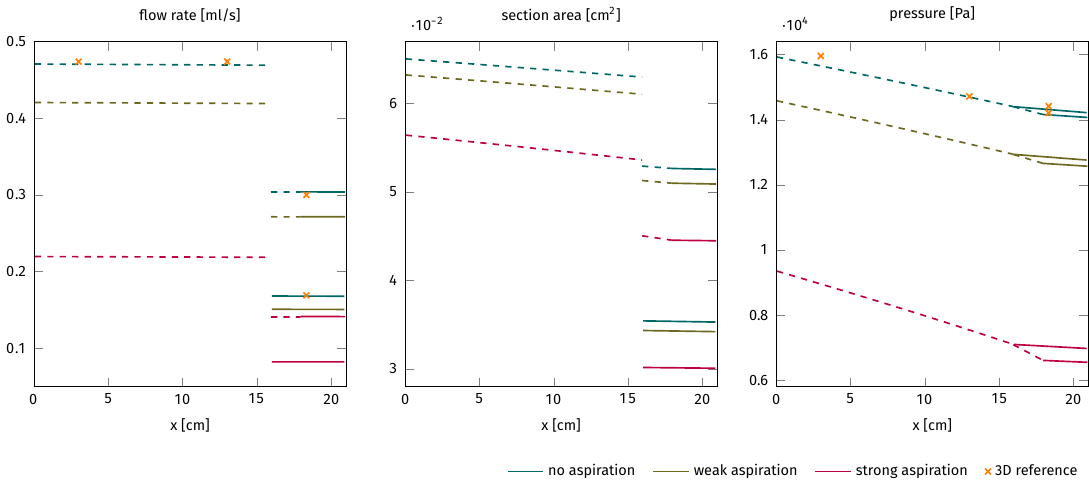}}
  \caption{Model results for the test case in Figure~\ref{fig:bifurcation_schematic}(B), where the catheter passes through the bifurcation point. Flow rate, section area, and pressure over the vessels are presented. In the catheterized segments (represented by dashed lines) total flow rate and section area are shown that include the flow and the size of the catheter. Either no aspiration, weak aspiration (suction with flow rate $-0.05\text{ml/sec}$) and strong aspiration (suction with flow rate $-0.25\text{ml/sec}$) is applied.}\label{fig:aspiration_bifurcation_2}
\end{figure*}

Simulation results of the setting in Figure~\ref{fig:bifurcation_schematic}B, where the catheter extends into \textit{vessel~2}, are presented in Figure~\ref{fig:aspiration_bifurcation_2}.
In this case a larger pressure drop in \textit{vessel~2} is predicted by the model, which along with the flow rate profiles agrees with the 3D and 1D literature results in the case that no aspiration was applied. Similarly to the previous case, the application of aspiration leads to a reduction of the flow rate in the arteries, as well as a reduction of the section area in the entire length of the catherized segments, which now extend to the daughter vessel. 
Aspiration also increases the relative pressure drop between the \textit{vessels~3} and \textit{2}. 

\subsection{Relation between suction force and wall distance}\label{sec:distancestudy}
In this section we apply the model to a blocked artery and particularly study how the distance of the catheter tip to the occlusion affects the exerted force due to aspiration thrombectomy. We adopt the setup from~\cite{shi2017suction} and compare simulation results by our reduced model with the 3D CFD computations the authors obtained. The suction force--suction distance relation was studied using a cylindrical vessel (inner diameter $D_0=0.4$ cm) representing the middle cerebral artery, with a concentrically aligned catheter (inner diameter $0.104$ cm, outer diameter $0.203$ cm) inserted from the proximal end to a depth of $5.7$ cm. A rigid thrombus, simply modeled as a domain boundary, was positioned in the vessel behind the catheter, and the suction distance $L$ (the gap between the catheter tip and the proximal thrombus surface) was varied from $0.01$ to $0.5$ cm. The proximal vessel inlet pressure was set to $12$ kPa and the catheter inlet suction pressure to $-30$ kPa. Moreover, blood was modeled with density $\rho=1.3$ g/ml and dynamic viscosity $\mu=0.032$ Pa\,s. The suction force S is defined as the sum of all pressure forces acting on the proximal and distal surfaces of the thrombus, computed as
\begin{equation}\label{eq:suctionforce}
  S = F_{pp} - F_{pd} = p_{pp} \,A_{pp} - p_{pd} \, A_{pd},
\end{equation}
where $p_{pp}$ and $A_{pp}$ are the proximal-side pressure and section area obtained from the simulation and $p_{pd}=12$ kPa and $A_{pd}= D_0^2 \frac{\pi}{4}$ are the fixed distal pressure and section area. In the special case of full contact ($L = 0$), the force is calculated as $S = F_{pv} + F_{pc} - F_{pd}$, accounting separately for the vascular pressure on the annular area and the catheter suction pressure on the inner circular area.

\begin{figure*}
  \iftoggle{usetikz}{\input{occlusion-aspiration.tikz}}{\includegraphics{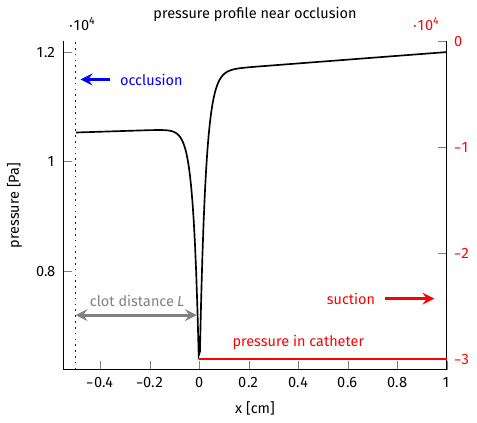}}
  \begin{minipage}[b]{0.52\linewidth}%
    \centering
    \includegraphics[trim={1cm 0 0 0},clip, width=\linewidth]{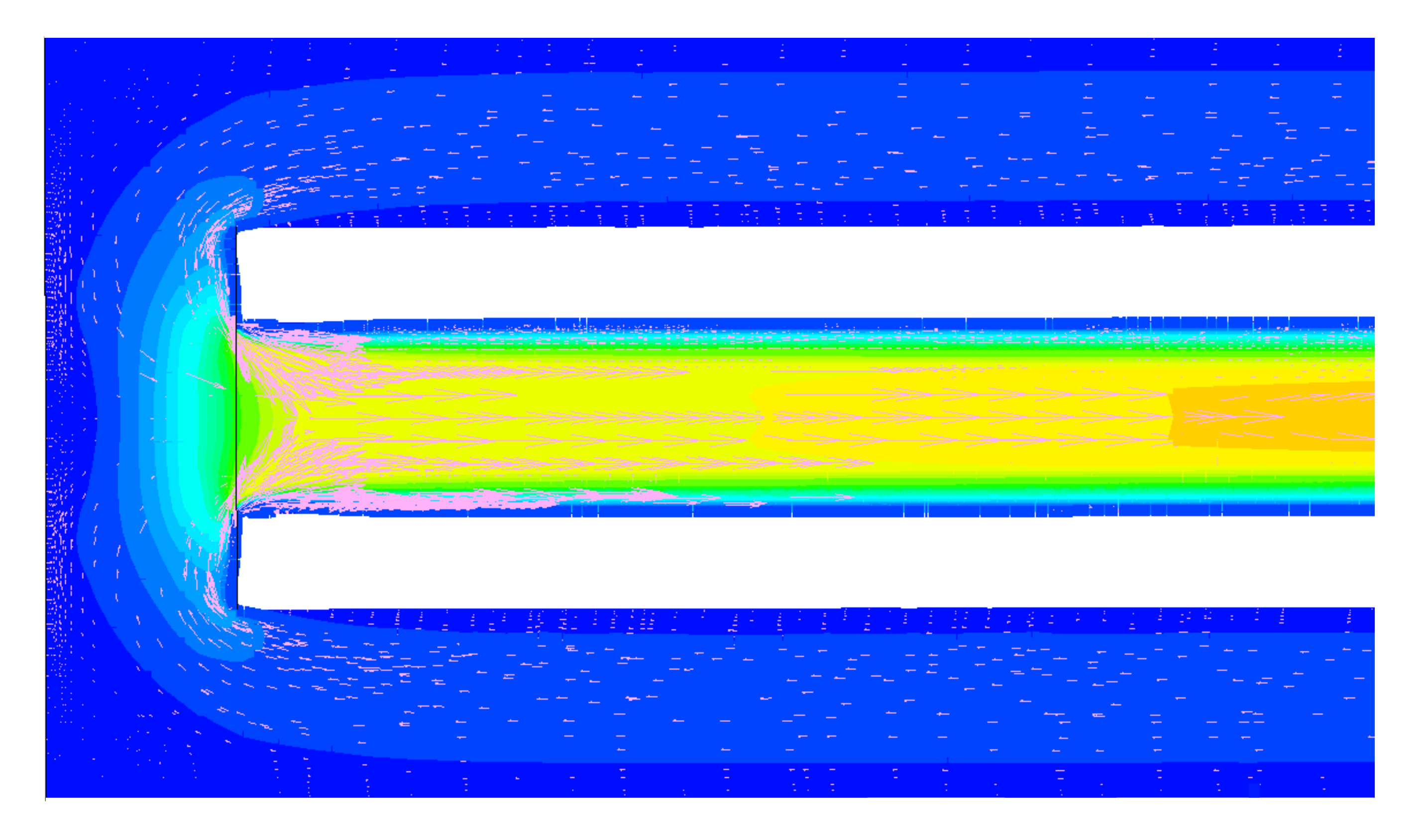}\\
    \iftoggle{usetikz}{\input{colorbar.tikz}}{\includegraphics{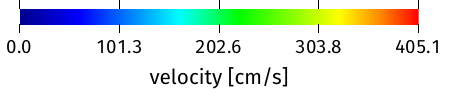}}
  \end{minipage}
\caption{Simulated flow profiles for aspiration thrombectomy with catheter tip close to an occlusion. Left: Pressure over the partially catheterized segment computed by the 1D model with indicated clot position, clot distance and negative pressure in the catheter. At the catheter tip the pressure attains a sharp minimum that smoothly connects with the pressure level in the catheterized and the reduced pressure level in the catheter-free segment near the occlusion. Right: Steady velocity profile with modeled suction through a catheter from a 3D CFD simulation.}\label{fig:occlusion-aspiration}
\end{figure*}

To implement the blocked artery scenario in the 1D model a partially catheterized configuration in analogy to Section~\ref{sec:partiallycatheterized} is considered, with the difference that at the outlet a resistive boundary condition with a reflection coefficient of $R_T=0.99$ is imposed to account for the occlusion. For consistency between the employed models we set the inlet to the right boundary and the outlet to the left boundary; the catheter tip is positioned at $x=0$ cm. Unlike the 3D CFD simulation, in which vessel and catheter are rigid and thus $A_{pp} = D_0^2 \frac{\pi}{4}$ in~\eqref{eq:suctionforce} is constant, the 1D model assumes an elastic vessel and catheter. Young's modulus of the vessel and the catheter and the external pressure are chosen similarly as in Section~\ref{sec:partiallycatheterized} and for efficiency of the simulations we employ the catheter equilibrium approximation. As the experiment addresses dynamics on a small spatial scale visco-elasticity of the vessel is assumed to plays a role, which is why we employ the splitting approach from Section~\ref{sec:splitting} with a positive dynamic viscosity $\gamma_c$. Some of the parameters have been fitted using an meta-heuristic algorithm, for details see Appendix~\ref{sec:numsetup}.

Independent of the clot distance $L$ the transient 1D simulations reach a steady state after a period of $200$ ms. Figure~\ref{fig:occlusion-aspiration} (left) exhibits the steady pressure profile in case of $L=0.5$ cm. While the pressure in the catheter is constant, the pressure in the vessel attains a sharp minimum at the catheter tip that smoothly connects with the pressure levels in the catheterized and the catheter-free segment near the occlusion. As $L$ is reduced the pressure at the occlusion approaches this minimum, hence, leading to a reduction in suction force. Figure~\ref{fig:occlusion-aspiration} (right) shows the velocity profile from an analogue 3D CFD simulation.

\begin{figure*}
  \iftoggle{usetikz}{\input{distance-relation.tikz}}{\includegraphics{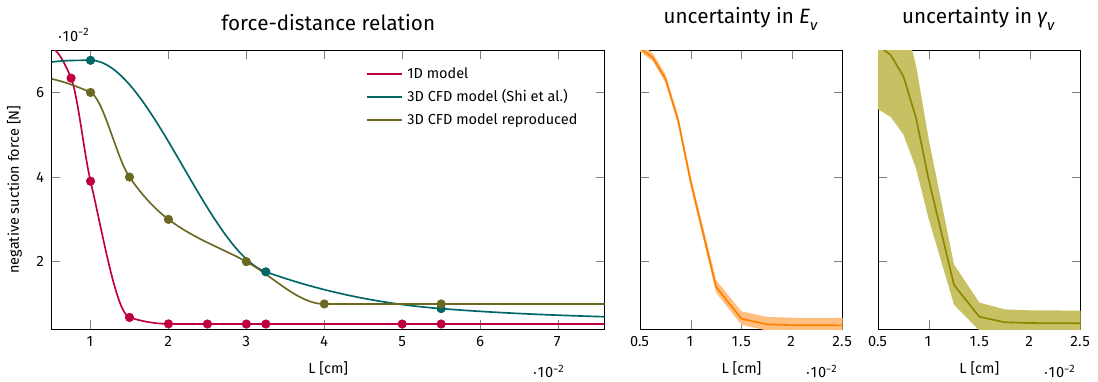}}
  \caption{Relation between the suction force $S$ and the distance $L$ between the catheter tip and the occlusion. Left: Comparison between the results of the reduced 1D model and 3D CFD simulations. Data points are interpolated by cubic splines. As data has been sparse in the critical regime of small $L$ in \cite{shi2017suction} complementary 3D CFD simulations were conducted. The 1D model with elastic walls qualitatively reproduces the force-distance relation in the 3D model with rigid walls. Right: The effect of uncertainty in Young's modulus $E$ and the dynamic wall viscosity $\gamma_v$ on the force distance relation represented by interquartile uncertainty bands. While variation in $E$ mainly affects the suction force plateau for larger $L$, variation in $\gamma_v$ might also lead to changes of the suction force for small $L$ and the decrease of characteristics of the force-distance relation as $L$ increases.}\label{fig:distance-relation}
\end{figure*}

For large $L$ the 1D model reproduces the 3D CFD predictions of the suction force with reasonable accuracy; for any $L\geq 0.76$~cm the relative difference in suction force between the interpolated 1D model and literature results is below $25$ \%. Since in literature only few data points are available for smaller $L$ we have reproduced the 3D CFD results and, in particular, computed further predictions for smaller catheter tip distances, see Appendix~\ref{sec:cfd} for details. Figure~\ref{fig:distance-relation} (left) compares the force--distance relations obtained by the 1D model to the 3D CFD predictions for small $L$. Also, in this regime the 1D model qualitatively reproduces the force-distance relation of the 3D model. While a difference with respect to the initial decrease of the force is given between the models, a discrepancy is expected due to the different model assumptions on elasticity of both the vessel and the catheter. 

Major benefit of the 1D model is the reduction of computation time over high-fidelity 3D simulations. We leverage this property and use Monte-Carlo simulations, see Section~\ref{sec:montecarlo}, to analyze the effect of uncertainty in the estimated model variables on the force-distance relation. In particular, we assigned the estimated Young's modulus $E_v$ and the dynamic wall viscosity $\gamma_v$ an independent uniform distribution over a relative range of $30\%$ and $40\%$, i.e., over the intervals $[0.7 E, 1.3 E]$ and $[0.6 \gamma_v, 1.4 \gamma_v]$, respectively. Based on $500$ samples and a discretization of the suction distance over $15$ points we computed the mean force--distance relation along with the central 50\% uncertainty band for both parameters, which are presented in Figure~\ref{fig:distance-relation} (right two panels). The results show that firstly variation in $E_v$ mainly affects the suction force plateau for larger $L$, and secondly variation in $\gamma_v$ might also lead to changes of the suction force for small $L$ and the decrease characteristics of the force-distance relation as $L$ increases.

\section{Discussion}
In this work, we have presented a reduced one-dimensional model for aspiration in blood vessels that accounts for the elastic behavior of both the vessel wall and the catheter. The aim was to provide a computational framework capable of capturing the essential hemodynamic phenomena during aspiration thrombectomy such as, suction-driven pressure drops and the coupling between catheterized and uncatheterized vessel segments—without the computational burden associated with full three-dimensional CFD simulations. The proposed model enables numerical computations that are substantially more efficient than high-fidelity 3D approaches, with the potential to simulate faster than real time, making it a viable tool for preoperative planning. A further practical advantage of the 1D formulation is the simplicity of geometry preparation: patient-specific anatomies can be discretized from minimal data, requiring only vessel connectivity, centerline lengths, and local diameters, thereby avoiding the labor-intensive mesh generation pipeline inherent to 3D simulations. The present work extends the earlier contribution \cite{pradhan2024} on 1D blood flow in a vessel with a rigid catheter, in which aspiration thrombectomy was identified as a promising future application; here, we realize that extension by introducing vessel and catheter elasticity, formulating the segment-to-segment coupling conditions, and benchmarking the model against both 3D CFD results and reference data from the literature.

Two modeling features introduced in this work merit particular discussion. The first is the inclusion of vessel wall viscoelasticity, which transforms the governing equation for the flow rate from a purely hyperbolic into a parabolic form. Our results demonstrate that this distinction becomes critical when the dynamics of interest occur on the millimeter scale, as is the case for the localized pressure field near the catheter tip. The suction-force--suction-distance analysis in Section~\ref{sec:distancestudy} showed that the viscoelastic model reproduces the sharp, diffuse pressure drop profile predicted by 3D CFD simulations, whereas a purely hyperbolic formulation, while still adequate for resolving dynamics on the centimeter scale, is systematically incapable of reproducing this near-tip pressure gradient. Importantly, this improved accuracy does not come at a prohibitive computational price: the implicit-explicit splitting strategy proposed in Section~\ref{sec:splitting} allows the viscoelastic terms to be treated implicitly while the hyperbolic terms remain explicit, so that the overall time-step size is governed by the latter. As a result, the viscoelastic model required only negligible additional cost relative to the purely hyperbolic formulation. The second feature is the simplified catheter model, which replaces the full elastic dynamics of the catheter wall with an equilibrium state determined by the model parameters and boundary conditions. The comparison in Section~\ref{sec:partiallycatheterized} confirmed that this catheter equilibrium approximation reproduces the fully elastic catheter model with high accuracy for time scales exceeding a few microseconds, provided that the Young's modulus of the catheter material is sufficiently large-a condition that is satisfied for typical clinical catheters. The practical benefit of this simplification is substantial: in our simulations it yielded a reduction in computation time by factors up to $65$, confirming the suitability of the overall framework for time-critical applications such as treatment planning.

A further advantage of the present framework stems from the numerical relaxation approach underlying the Lax--Friedrichs-type finite volume scheme described in Section~\ref{sec:spacedisc}. Originally introduced in~\cite{herty2023centrschemtwo} and subsequently adopted and verified for 1D hemodynamic modeling in \cite{beckers2025laxfried}, this relaxation serves two purposes: it yields a robust and simple finite volume discretization of the hyperbolic subsystem, and, more importantly, it considerably simplifies the coupling of different model segments within a vascular network. In the present work, this property has been exploited to construct straightforward nodal solvers at two distinct locations: catheterized bifurcations, where three or more vessel segments meet and the catheter may span one or more of the outgoing branches, and the catheter end, where the transition between the catheterized and uncatheterized lumen occurs. For the latter, we derived closed, purely algebraic expressions for the coupling data, eliminating the need for iterative solvers at this interface. The resulting coupling is highly efficient and, as demonstrated by the validation cases, produces simulation results free of spurious boundary layers at the segment interfaces. We acknowledge, however, that the current framework treats vascular occlusions only in a static sense, representing them as domain boundaries rather than as dynamic entities accounted for in the model equations. A more comprehensive model would account for mobile clots that may be partially or fully aspirated into the catheter lumen, fragment during the procedure, or embolize to downstream vessel segments. Such clot dynamics would necessitate additional transport and constitutive modeling but can be accommodated within the existing coupling architecture, as the nodal solvers and segment-based formulation are sufficiently general to incorporate extended hyperbolic systems. The development of these enhanced occlusion models is a natural direction for future work.

While the comparisons with data from the literature and with established 3D CFD simulations presented in this work provide strong evidence for the accuracy of the proposed model, they do not constitute a complete validation. Ideally, future studies will include in vitro experiments, for instance, flow phantoms with controlled aspiration conditions, and, where feasible, in vivo data from clinical aspiration thrombectomy procedures. It is worth noting that this limitation is not specific to the present 1D framework; even state-of-the-art 3D CFD and FSI models of aspiration thrombectomy still lack comprehensive clinical validation. Nevertheless, we believe that the proposed approach is promising for clinical application, particularly in supporting patient-specific decisions such as optimal cannula tip positioning during the intervention. The computational efficiency is a key feature in this context: as demonstrated by the uncertainty quantification study in Section~\ref{sec:distancestudy}, the model makes it feasible to perform large numbers of simulations required to investigate the sensitivity of clinically relevant quantities to variations in anatomical and procedural parameters.
With further refinement, particularly the incorporation of dynamic clot models and validation against experimental data, we believe that reduced 1D models can become a practical computational tool for aspiration thrombectomy, complementing high-fidelity simulations and ultimately contributing to improved procedural outcomes.

\section{Conclusion}
We have presented a reduced one-dimensional model for aspiration in blood vessels that accounts for the elasticity of both the vessel wall and the catheter. The inclusion of vessel wall visco-elasticity yields a parabolic formulation for the flow rate that accurately captures the localized pressure dynamics at the catheter tip on the millimeter scale, as confirmed by comparison with 3D CFD simulations, while a splitting strategy ensures only negligible additional computational cost over the hyperbolic formulation. A simplified catheter equilibrium approximation was shown to reproduce the fully elastic catheter model with high accuracy while reducing computation time significantly. The numerical relaxation approach used in the model enables a robust Lax--Friedrichs-type finite volume scheme and, through closed algebraic nodal solvers, straightforward coupling between catheterized and uncatheterized vessel segments including bifurcations and the catheter tip. The model was benchmarked against literature data and 3D CFD results in multiple test cases. Its computational efficiency makes it suitable for time-critical applications such as uncertainty quantification. Future work will focus on extending the framework to incorporate dynamic clot models during aspiration.

\subsection*{Funding}
The authors thank the Deutsche Forschungsgemeinschaft (DFG, German Research
Foundation) for the financial support through SPP
2311 (Robust Coupling of Continuum-Biomechanical In Silico Models to Establish Active
Biological System Models for Later Use in Clinical Applications -- Co-Design of Modelling,
Numerics and Usability) within the Project 548864771, through SPP 2410 (Hyperbolic Balance Laws in Fluid
Mechanics: Complexity, Scales, Randomness) within the Project 525842915 and under 320021702/GRK2326 (Graduate College Energy, Entropy, and Dissipative Dynamics).

\subsection*{Ethics Statement}
The authors have nothing to report.

\subsection*{Conflicts of Interest}
The authors declare no conflicts of interest.

\subsection*{Data Availability Statement}
The data that support the findings of this study are available from the corresponding author upon reasonable request.

\subsection*{Declaration on the Use of Generative AI Tools}
During the preparation of this work, the authors used the generative AI tool GLM 5.3 (Zhipu AI) solely for language editing in the Introduction, Discussion and Conclusion sections. The tool was not used for the development of code, the execution of simulations, the analysis of results, or the generation of any scientific content. All AI-assisted text was reviewed and edited by the authors, who take full responsibility for the integrity and accuracy of the published work.

\appendix
\setcounter{section}{0}
\renewcommand{\thesection}{A\arabic{section}}

\section{Solver for the hyperbolic coupling problem}\label{sec:hcoupling}
We consider Jin-Xin relaxation of the systems \eqref{eq:conservativecatheter} and \eqref{eq:conservativeAu} with relaxation speed $\lambda>0$ regarding the quantities $A_n$, $u_n$, $A_v$ and $u_v$. The catheter related quantities $A_c$ and $u_c$ are relaxed with the speed $\lambda_c$ that is generally significantly larger than $\lambda$. For the additional relaxation variables we impose the new coupling conditions
\begin{align}
  v^{A_n}(0^-, t) + v^{A_c}(0^-, t) &= v^{A_v}(0^+, t), \label{eq:couplingvA}\\
  v^{u_n}(0^-, t) + \frac{u_v(0^+, t)^2 - u_n(0^-, t)^2}{2} &= v^{u_c}(0^+, t) \label{eq:couplingvu}
\end{align}
that are consistent with the original conditions \eqref{eq:couplingA} and \eqref{eq:couplingu}. In more details, if the relaxation variables are replaced by the corresponding fluxes in the conservative systems, these conditions are equivalent to the original ones.
  
In the following we design a Riemann solver
assigning the trace data $A_n^-$, $u_n^-$, $A_c^-$, $u_c^-$, $v^{A_n-}$, $v^{u_n-}$, $v^{A_c-}$, $v^{u_c-}$, $A_v^+$, $u_v^+$, $v^{A_v+}$, $v^{u_v+}$ to the coupling data $[A_n]_R$, $[u_n]_R$, $[A_c]_R$, $[u_c]_R$, $v^{A_n}_R$, $v^{u_n}_R$, $v^{A_c}_R$, $v^{u_c}_R$, $[A_v]_L$, $[u_v]_L$, $v^{A_v}_L$, $v^{u_v}_L$ at the interface. The coupling data both solve the half-Riemann problems with respect to the two relaxation systems at the interface and satisfy the conditions \eqref{eq:couplingA}, \eqref{eq:couplingu}, \eqref{eq:couplinghomogeneous}, \eqref{eq:couplingvA} and \eqref{eq:couplingvu}.

From the homogeneous conditions \eqref{eq:couplinghomogeneous} for the device related variables we obtain the coupling states
\begin{equation*}
  [A_c]_R = A_c^-, \quad [u_c]_R = u_c^-, \quad v^{A_c}_R = v^{A_c-}, \quad v^{u_c}_R = v^{u_c-}.
\end{equation*}
To identify the remaining coupled section area variables we note that the Lax curves of the relaxation systems imply the relations
\[
  v^{A_n}_R = v^{A_n-} + \lambda (A_n^- -[A_n]_R), \quad  v^{A_v}_L = v^{A_v+} + \lambda ([A_v]_L - A_v^+) .
\]
Thus~\eqref{eq:couplingvA} allows us to deduce a linear relation between the unknowns $ [A_n]_R$ and $[A_v]_L$, from which we obtain the states
\begin{align*}
    2\,[A_v]_L &= A_n^-  + A_v^+ + \frac{1}{\lambda} \left(v^{A_n-} + v^{A_c}_R - v^{A_v+}\right) \\
                      &\quad + \left(\sqrt{\frac{A_c}{\pi}} + h_c \right)^2 \pi,\\
  [A_n]_R &= A_n^- + A_v^+ - [A_v]_L \\
  &\quad + \frac{1}{\lambda} \left(v^{A_n-} + v^{A_c}_R - v^{A_v+}\right).
\end{align*}
Assuming that the trace data satisfies the coupling conditions, particularly \eqref{eq:couplingvA}, up to a small error we can infer that $v^{A_n-} + v^{A_c}_R - v^{A_v+}$ is small and hence the coupling states are positive and physically admissible. 

Next, we determine the velocities  $[u_n]_R$ and $[u_v]_L$ and their corresponding relaxation states $v^{u_n}_R$ and $v^{u_c}_R$. We parameterize these states on the Lax-curves of the relaxation systems as 
\begin{align*}
  [u_n]_R &= u_n^- - \sigma_u^-, \quad &&v_R^{u_n} = v^{u_n-} + \lambda \sigma_u^-,\\
  [u_v]_L&= u_v^+ + \sigma_u^+, \quad &&v_L^{u_v} = v^{u_v+} + \lambda \sigma_u^+.
\end{align*}
Substituting this parametrization in~\eqref{eq:couplingu} we obtain
\begin{equation}\label{eq:linearsigmau}
  \begin{aligned}
  \sigma_u^- &= \frac{s_0}{[A_n]_R} - \frac{[A_v]_L}{[A_n]_R} \sigma_u^+, \\
  s_0 &\coloneqq [A_n]_R \, u_n^- + [A_c]_R \, [u_c]_R - [A_v]_L \, [u_v^+] .
  \end{aligned}
\end{equation}
Next, employing the parametrization in \eqref{eq:couplingvu} yields
\begin{align*}
 \lambda (\sigma_u^+ - \sigma_u^-) - u_v^+ \sigma_u^+  - u_n^- \sigma_u^- + \frac 1 2 ({\sigma_u^-}^2 - {\sigma_u^+}^2) \\= v^{u_n-} + \frac 1 2 ({u_v^{+}}^2 - {u_n^-}^2) - v^{u_v+}.
\end{align*}
Substituting now \eqref{eq:linearsigmau} and setting $q= \frac{[A_v]_L}{[A_n]_R}$ we derive the quadratic equation
  \begin{equation}\label{eq:sigmaquadratic}
  \begin{split}
    \frac 1 2  \left[ q^2 - 1\right] {\sigma_u^+}^2 + \left[ \lambda \left( 1 + q \right) + u_n^- q - u_v^+ - \frac{s_0 q}{[A_n]_R}\right] \sigma_u^+ \\
    + \left(\frac{s_0}{2[A_n]_R} - \lambda - u_0^-\right) \frac{s_0}{[A_n]_R} - s_1 = 0.
  \end{split}
\end{equation}
Due to \eqref{eq:couplingA} it holds that $[A_v]_L/[A_n]_R>1$ and therefore the quadratic term does not vanish. Moreover, we note that in practice the trace states $u_n^{-}$, $u_v^{+}$,  $v^{u_n-}$ and $v^{u_v+}$ are close to the corresponding coupling states and thus due to \eqref{eq:couplingu} and \eqref{eq:couplingvu} the terms $|s_0|$ and $|s_1|$ are relatively small. Consequently the constant term in \eqref{eq:sigmaquadratic} has a relative small absolute value in comparison to the linear term. For these reasons we can expect that~\eqref{eq:sigmaquadratic} has two real roots, from which we select the one with the smaller absolute value.

\section{Details on the numerical simulations}\label{sec:numsetup}
In this appendix, we provide additional details on the numerical simulations of the 1D model. Those were conducted either with the first order scheme~\eqref{eq:lxf} or with a second order extension based on the MUSCL scheme and the Minmod-limiter, see~\cite{beckers2025laxfried, leer1979towarultimconserdifferschem} for details. The time increment was updated in every time step with Courant number CFL$=0.9$ for the first order and CFL$=0.24$ for the second order scheme. The second order scheme was used to produce the numerical results in Section~\ref{sec:partiallycatheterized} and the first order scheme in the other parts of the paper.

In Section~\ref{sec:partiallycatheterized} we chose $\Delta x=0.01$ cm except for the convergence study in Figure~\ref{fig:elasticity-study} (right), where the spatial resolution was reduced to $\Delta x=0.05$ cm. We note in this respect, that the fully elastic model took between $4$ ($E_c=1$ MPa) to $65$ ($E_c=500$~MPa) times the computation time of the catheter equilibrium approximation. 

In Section~\ref{sec:bifurcationstudy} we choose different spatial discretizations for each vessel segment, in each one of the two configurations shown in Figure~\ref{fig:bifurcation_schematic}. For the first case shown in Figure~\ref{fig:bifurcation_schematic}A, $\Delta x= 0.02, 0.125$ and $0.125$ cm for \textit{vessel~1, 2 \& 3}, respectively, and for the second case shown in Figure~\ref{fig:bifurcation_schematic}B, $\Delta x= 0.1067, 0.01$ and $0.125$ cm for \textit{vessel~1, 2 \& 3}, respectively. In the numerical simulations of the second scenario, we split \textit{vessel~2} into two segemnets, a catheterized and a non-catheterized one, to account for the catheter tip position, and apply the coupling conditions at the interface between the two segments, as described in Sections~\ref{sec:coupling-problem} and \ref{sec:coupling-problem-catheterized}. 
The equlibrium state of the system is reached at time $T=5$ s. 

In Section~\ref{sec:distancestudy} the computational domain of the catheterized part of the segment was reduced to $1$ cm, which did not affect the simulation results but allowed for more efficient computations. We have employed a fine grid with $\Delta x=0.00125$ cm. Young's modulus of the vessel and the catheter were chosen $E_v =641.484$ kPa and $E_c=30$ MPa. Moreover, we set $\gamma_c = 8.97105$ Pa\,s and chose the external pressure $p_v^\text{ext}= p_n^\text{ext}=4.3$~kPa and $p_c^\text{ext}=490$~kPa. These parameters were partially tuned using the adaptive differential evolution optimizer from \cite{storn1997differ} to fit the literature results for fixed $L=0.5$ cm.

\section{Details on the comparative CFD simulations}\label{sec:cfd}
A 3D CFD model was created with the same setup settings as presented in~\cite{shi2017suction}. In brief, a steady state CFD analysis was performed with an incompressible Newtonian fluid with a density $\rho=1.3 $g/ml and a dynamic viscosity $\mu= 3.2$ mPa\,s. The circular vessel had an inner diameter of $0.4$ cm and a length of 6 cm (see Figure~\ref{fig:catheterized}). The catheter had an inner diameter of $0.104$ cm and an outer diameter of $0.203$ cm. A negative suction pressure of $-30$ kPa was imposed at the catheter outlet and a positive pressure of $12$ kPa was imposed at the vessel inlet. Simulations were performed at $6$ different tip distances, $0.01$ cm, $0.015$ cm, $0.02$ cm, $0.03$ cm, $0.04$ cm and $0.055$ mm. The suction force acting on the thrombus was determined from the difference of area averaged pressure acting on the thrombus and an assumed distal thrombus pressure of $12$ kPa. An unstructured tetrahedral mesh was used with approximately $2 \times 10^5$ elements and laminar flow conditions were assumed. All simulations were performed with ANSYS Fluent 2024R1\footnote{Ansys Inc., Canonsburg, PA} using the coupled pressure-velocity coupling scheme and a convergence target of $10^{-3}$ in terms of root mean square for continuity and momentum equations. 

\bibliographystyle{abbrvurl} 
\bibliography{aspiration-refs}
\end{document}